\documentclass[12pt]{article}
\usepackage{amssymb,amsmath,amsthm,amsopn}
\numberwithin{equation}{section}

\newtheorem{theorem}{Theorem}
\newtheorem{proposition}{Proposition}[section] 
\newtheorem{corollary}[proposition]{Corollary}
\newtheorem{lemma}[proposition]{Lemma}
\newtheorem{definition}[proposition]{Definition}
\newtheorem{example}[proposition]{Example}

\let\ds=\displaystyle
\def\coun{\varepsilon}
\def\bn{\begin{equation}}
\def\ed{\end{equation}}

\def\<{\langle}
\def\>{\rangle}
\newcommand{\CC}{{\mathbb C}}

\newcommand{\ZZ}{{\mathbb Z}}

\def\r#1{(\ref{#1})}

\def\ot{\otimes}

\newcommand{\Uqdva}{{U_{q}^{}(\widehat{\mathfrak{sl}}_{2})}}
\def\sk#1{\left(#1\right)}
\def\skk#1{\left[#1\right]}

\def\Pep{{P}^+_e}
\def\Pem{{P}^-_e}
\def\Pepm{{P}^\pm_e}

\def\Pfpm{{P}^\pm_f}
\def\Pfp{{P}^+_f}
\def\Pfm{{P}^-_f}

\def\Uqgln{U_q(\widehat{\mathfrak{gl}}_N)}
\def\Uqglnd{\overline{U}_q(\widehat{\mathfrak{gl}}_N)}

\def\Uqd#1{\overline{U}_q(\widehat{\mathfrak{gl}}_{#1})}

\def\Uqsl2{U_q(\widehat{\mathfrak{sl}}_2)}

\def\Uqgldva{\overline{U}_q(\widehat{\mathfrak{gl}}_2)}
\def\Uqgltwo{{U}_q(\widehat{\mathfrak{gl}}_2)}
\def\Uqbm{U_q(\mathfrak{b}^-)}
\def\Uqbp{U_q(\mathfrak{b}^+)}
\def\ff{{F}}

\def\EE{{\rm E}}\def\FF{{\rm F}}

\def\Uqbp{U_q(\mathfrak{b}^+)}

\def\si{\sigma}
\def\RR{{\rm R}}
\def\LL{{\rm L}}

\def\E{{\sf E}}
\def\F{{\mathcal{F}}}

\def\qsym{\varpi}
\def\tSym{\overline{\rm Sym}}

\def\lll{l}
\def\rr{r}
\def\ss{s}
\def\mm{m}

\def\seg#1#2{[#2,#1]}
\def\segg#1{[#1]}
\def\meg#1#2{[#1,#2]}
\def\megg#1{[#1]}

\def\Ff{f}

\let\qsym=\beta

\def\hS{\hat{S}}

\def\qed{\hfill$\square$\medskip}
\def\tS{\hat{S}}
\def\emb{\Psi}

\def\Ags{\mathcal{E}}
\def\Bgs{\mathcal{D}}

\def\tSgs{\tilde{\mathcal{S}}}
\def\simm{\thickapprox}
\def\simJ{\sim_{\!\scriptscriptstyle{J}}}
\def\simI{\sim_{\!\scriptscriptstyle{I}}}
\def\simK{\sim_{\!\scriptscriptstyle{K}}}
\def\simJK{\sim_{\!\scriptscriptstyle{J,K}}}
\def\cU{{U}}
\def\fb{{\mathfrak b}}
\def\Zfun{\mathbb{Z}}
\def\Yfun{\mathbb{Y}}
\def\Rfun{\mathbb{R}}
\def\Aint{\mathbb{A}}
\def\Qfun{\mathbb{Q}}

\begin{document}

\begin{center}

\hfill ITEP-TH-40/08\\
\hfill LAPTH-1270/08\\
\bigskip
{\Large\bf Bethe Ansatz for the Universal Weight Function}
\par\bigskip\medskip
{\bf
Luc Frappat$^\circ$\footnote{E-mail: frappat@lapp.in2p3.fr},
Sergey Khoroshkin$^{\star}$\footnote{E-mail: khor@itep.ru},
Stanislav Pakuliak$^{\star\bullet}$\footnote{E-mail: pakuliak@theor.jinr.ru},
\'Eric Ragoucy$^\circ$\footnote{E-mail: ragoucy@lapp.in2p3.fr}}
\par\bigskip\medskip
$^\star${\it Institute of Theoretical \& Experimental Physics, 117259
Moscow, Russia}
\par\smallskip
$^\bullet${\it Laboratory of Theoretical Physics, JINR,
141980 Dubna, Moscow reg., Russia}
\par\smallskip
$^\circ${\it Laboratoire d'Annecy-le-Vieux de Physique Th\'eorique
LAPTH, \\
UMR 5108 du CNRS associ\'ee \`a l'Universit\'e de Savoie\\
BP 110, F-74941, Annecy-le-Vieux cedex, France}\\
\bigskip
\end{center}

\thispagestyle{empty}

\begin{abstract}
We consider universal off-shell Bethe vectors given in terms of Drinfeld
realization of the algebra  $\Uqglnd$ \cite{KPT,KhP-GLN}. We investigate
ordering properties of the product of the transfer matrix and these vectors.
  We derive that these vectors are eigenvectors of the transfer matrix
  if their Bethe parameters satisfy the universal Bethe equations
  \cite{ACDFR}.
\end{abstract}


\setcounter{footnote}{0}

\section{Introduction}

Algebraic Bethe ansatz for  quantum integrable models with $\mathfrak{gl}_2$
symmetry as well as hierarchical Bethe ansatz for models with higher rank
symmetries solves the eigenvalue problem for the set of the commuting quantum
integrals of motion. The eigenvectors in these methods are built from the
matrix elements of the monodromy matrix which satisfies Yang-Baxter relation
defined by the quantum $\RR$-matrix. Quantum integrable models solvable by
these methods correspond to the different monodromy matrices and quantum
$\RR$-matrices. Monodromy matrices can be obtained by considering the different
representations of the quadratic algebras which have the same type of defining
relations as monodromy matrices have. In this case one can use the generating
series of the elements of this quadratic algebra as kind of the {\it universal}
monodromy matrix and try to reformulate the Bethe ansatz on the universal level
without specification to any concrete representations or concrete integrable
model.

Such an approach to find the eigenvalues for the quantum integrable models with different
boundary conditions and symmetries was elaborated in \cite{ACDFR} using certain
analytical assumptions on the structure of these eigenvalues. This method
is called
analytical Bethe ansatz and by construction was unable to yield an information on the structure
of the corresponding eigenvectors. In case of the models with $\mathfrak{gl}_N$ symmetry
the method to build the eigenvectors from the matrix elements of the monodromy matrices
was designed in \cite{VT,VT1,MTV,VT2} generalizing an approach of the hierarchical Bethe
ansatz formulated in \cite{KR83}. Authors of these papers used the universal monodromy which
satisfies the commutation relations of the $\mathfrak{gl}_N$ Yangian or Borel subalgebra
in the quantum affine algebra $\Uqgln$ to construct the universal Bethe vectors of the
hierarchical Bethe ansatz in terms of the matrix elements of the corresponding fundamental
$\LL$-operators.

Quantum affine algebra and their rational analogs -- the Yangian doubles have an alternative
to the $\LL$-operator description \cite{RS} in term of currents \cite{D88,DF}. For the quantum
affine algebras, it was proved
in \cite{KT} that the modes of Drinfeld currents coincide with Cartan-Weyl
generators of these algebras constructed from the finite set of Chevalley generators.
One may address the question whether it is possible to construct the universal Bethe vectors
from the current generators of the quantum affine algebra which serves as symmetry for some
quantum integrable models. This problem was investigated in \cite{EKhP} on a
rather general level and it was shown that in order to build the universal Bethe vectors from
Cartan-Weyl or current generators of the algebra one has to consider different
types of Borel subalgebras in the quantum affine algebras.
In this paper it was suggested to construct universal off-shell Bethe vectors for arbitrary
quantum affine algebra as
certain projections of products of the currents onto intersections of Borel
subalgebras of different types.
The generating parameters of the currents become after this identification
the Bethe parameters. The papers \cite{KhP-GLN,OPS} contain detailed analysis of these projections
for quantum affine algebra $\Uqgln$. In particular, these projections are explicitly expressed via
entries of the fundamental monodromy matrix and are identified with off-shell Bethe vectors of the
nested Bethe ansatz \cite{KR83}.

An algebraic Bethe ansatz always uses a special vector which is annihilated by some ideal
in the symmetry algebra (bare vacuum) and  Bethe vectors are obtained by the application of the
universal Bethe vectors to this vector. From the representation theory point of view we will call
such bare vacuum a weight singular vector. The Cartan-Weyl generators have a good property: their
products can be ordered in a natural way. If we are able to express the commuting integrals as well as the
universal Bethe vectors in terms of these generators, we may rise the question: what is special
in the universal Bethe vectors if their Bethe parameters satisfy the universal Bethe equations.
In this paper we address this question for the quantum affine algebra $\Uqglnd$.
We found that the Cartan-Weyl ordering of the product of the universal transfer matrix and the
universal Bethe vector
produces the same Bethe vectors modulo elements of the ideal which annihilates the bare vacuum
if the Bethe parameters satisfy the universal Bethe equations from \cite{ACDFR}.

All our calculations are performed on the level of generating series and the main technical
trick which helps us to perform the ordering calculations is the identity \r{spelpr1}
which is a particular case of more general relations between off-shell Bethe vectors
obtained in the paper \cite{FKPR-GS} using the technique of the generating series.

The paper is composed as follows. In Section~\ref{sect2}, all necessary statements for the
different realizations of the quantum affine algebra $\Uqglnd$ are given
and the main assertion of the
paper is formulated as the Theorem~\ref{main-th}. Section~\ref{sect3}
collects propositions which describe the ordering of the generating series of the Cartan-Weyl or current
generators. There, an
identity \r{spelpr1} is formulated: it is a new type of hierarchical relation between universal off-shell
Bethe vectors expressed in terms of the current generators. Section~\ref{sect4} is devoted
to the inductive proof of the main Theorem~\ref{main-th}.

\section{Basic algebraic structures}\label{sect2}

\subsection{$\Uqglnd$ as a quantum double}

Let $q$ be  a complex parameter not equal to zero or to a
root of unity.
Let $\E_{ij}\in{\textrm{End}}(\CC^N)$ be a matrix with the only nonzero entry
equal to $1$
at the intersection of the $i$-th row and $j$-th column.
Let $\RR(u,v)\in{\textrm{End}}(\CC^N\ot\CC^N)\ot \CC[[{v}/{u}]]$,
\begin{equation}\label{UqglN-R}
\begin{split}
\RR(u,v)\ =\ &\ \sum_{1\leq i\leq N}\E_{ii}\ot \E_{ii}\ +\ \frac{u-v}{qu-q^{-1}v}
\sum_{1\leq i<j\leq N}(\E_{ii}\ot \E_{jj}+\E_{jj}\ot \E_{ii})
\\
+\ &\frac{q-q^{-1}}{qu-q^{-1}v}\sum_{1\leq i<j\leq N}
(u \E_{ij}\ot \E_{ji}+ v \E_{ji}\ot \E_{ij})
\end{split}
\end{equation}
be a trigonometric $\RR$-matrix associated with the vector
representation of ${\mathfrak{gl}}_N$.

We will consider an associative algebra $\Uqglnd$ with unit
as a quantum double \cite{D-qg} of its Borel subalgbera generated by the modes
$\LL^{+}_{i,j}[k]$, $k\geq 0$, $1\leq i,j\leq N$ of the
$\LL$-operators
$\LL^{+}(z)=\sum_{k=0}^\infty\sum_{i,j=1}^N \E_{ij}\otimes \LL^{+}_{i,j}[k]z^{-
k}$, $\LL^{+}_{j,i}[0]=0$, $1\leq i<j \leq N$
subject to the relations
\begin{equation}\label{L+}
\RR(u,v)\cdot (\LL^{+}(u)\ot \mathbf{1})\cdot (\mathbf{1}\ot \LL^{+}(v))=
(\mathbf{1}\ot \LL^{+}(v))\cdot (\LL^{+}(u)\ot \mathbf{1})\cdot \RR(u,v)
\end{equation}
with a standard coproduct
\begin{equation}\label{coL+}
\Delta \sk{\LL^{+}_{i,j}(u)}=\sum_{k=1}^N\ \LL^{+}_{k,j}(u)\otimes
\LL^{+}_{i,k}(u)\,.
\end{equation}
We denote this subalgebra
$U_q(\mathfrak{b}^+)\subset\Uqglnd$ and call it {\em a standard Borel subalgebra} of $\Uqglnd$.
In \r{L+}  $\mathbf{1}=\sum_{i=1}^N \E_{ii}$.

According to the general theory \cite{D-qg} the whole algebra $\Uqglnd$ is generated
by the modes of the $\LL$-operator $\LL^+(z)$ and by the modes of the dual $\LL$-operator
$\LL^{-}(z)=\sum_{k=0}^{-\infty}\sum_{i,j=1}^N \E_{ij}\otimes \LL^{-}_{i,j}[k]z^{-
k}$, $\LL^{-}_{i,j}[0]=0$, $1\leq i<j \leq N$. The dual Borel subalgebra $U_q(\mathfrak{b}^-)$
has the same
algebraic and coalgebraic properties \r{L+} and \r{coL+} with $\LL^+(z)$
replaced by $\LL^-(z)$ everywhere. The commutation relation between opposite Borel
 subalgebras
\begin{equation}\label{L+-}
\RR(u,v)\cdot (\LL^{+}(u)\ot \mathbf{1})\cdot (\mathbf{1}\ot \LL^{-}(v))=
(\mathbf{1}\ot \LL^{-}(v))\cdot (\LL^{+}(u)\ot \mathbf{1})\cdot \RR(u,v),
\end{equation}
can be calculated using the non-degenerated pairing between these subalgebras.

The standard description of the quantum affine algebra $\Uqgln$
with dropped gradation element and at vanishing central element  can be obtained
from the algebra $\Uqglnd$ by imposing one more relation
\begin{equation}\label{rest}
\LL^+_{i,i}[0]\LL^-_{i,i}[0]=1\qquad i=1,\ldots,N\,.
\end{equation}
Here we shall not assume this restriction. We shall require only invertibility of the zero modes
of the diagonal matrix elements of $\LL$-operators.

\subsection{Current realization of $\Uqglnd$}

To obtain the current realization of the algebra $\Uqglnd$
\cite{D88} one has to introduce, according to \cite{DF},  the Gauss coordinates of the
$\LL$-operators. There are two different ways of introducing Gauss decompositions
of $\LL$-operators. Each of these two possibilities leads to hierarchical
relations for the universal Bethe vectors given by the hierarchical Bethe ansatz.
 One type of hierarchy occurs when the smaller algebra $\LL$-operator is embedded in the
upper-left
corner of the bigger algebra $\LL$-operator. The second type corresponds to the embedding
in the down-right corner \cite{VT2,OPS}. In this paper we will use the
latter embedding. Then, Gauss coordinates
$\FF^{\pm}_{b,a}(t)$,
$\EE^{\pm}_{a,b}(t)$, $b>a$ and
$k^\pm_{c}(t)$ are given by the decomposition:
\begin{equation}\label{L-op}
\LL^\pm(z)=\sk{\mathbf{1}+\sum^N_{a<b}\FF^\pm_{b,a}(z)\E_{ab}}
\cdot\sk{\sum^N_{a=1}k^\pm_{a}(z)\E_{aa}}\cdot
\sk{\mathbf{1}+\sum^N_{a<b}\EE^\pm_{a,b}(z)\E_{ba}}
\end{equation}
that is to say
\begin{align}\label{GF1}
\LL^{\pm}_{a,b}(t)&=\FF^{\pm}_{b,a}(t)k^+_{b}(t)+\sum_{b<m\leq N}
\FF^{\pm}_{m,a}(t)k^+_{m}(t)\EE^{\pm}_{b,m}(t),\qquad a<b,\\
\label{GK1}
\LL^{\pm}_{b,b}(t)&=k^\pm_{b}(t) +\sum_{b<m\leq N} \FF^{\pm}_{m,b}(t)k^\pm_{m}(t)
\EE^{\pm}_{b,m}(t),\\
\label{GE1}
\LL^{\pm}_{a,b}(t)&=k^\pm_{a}(t)\EE^{\pm}_{b,a}(t)+\sum_{a<m\leq N}
\FF^{\pm}_{m,a}(t)k^\pm_{m}(t)\EE^{\pm}_{b,m}(t),\qquad   a>b\,.
\end{align}

Considering the linear combinations of the Gauss coordinates
\begin{equation}\label{DF2}
F_i(t)=\FF^{+}_{i+1,i}(t)-\FF^{-}_{i+1,i}(t)\,,\quad
E_i(t)=\EE^{+}_{i,i+1}(t)-\EE^{-}_{i,i+1}(t)
\end{equation}
and diagonal Gauss coordinates $k^\pm_i(t)$ one can obtain  the defining  relations \cite{D88,DF}:
\begin{gather}
(q^{-1}z-q^{}w)E_{i}(z)E_{i}(w)= E_{i}(w)E_{i}(z)(q^{}z-q^{-1}w)\, ,  \notag \\
(z-w)E_{i}(z)E_{i+1}(w)= E_{i+1}(w)E_{i}(z)(q^{-1}z-qw)\, ,  \notag \\
k_i^\pm(z)E_i(w)\left(k_i^\pm(z)\right)^{-1}=
\frac{z-w}{q^{-1}z-q^{}w}E_i(w)\, , \notag \\
k_{i+1}^\pm(z)E_i(w)\left(k_{i+1}^\pm(z)\right)^{-1}=
\frac{z-w}{q^{}z-q^{-1}w}E_i(w)\, , \notag \\
k_i^\pm(z)E_j(w)\left(k_i^\pm(z)\right)^{-1}=E_j(w),
\qquad {\rm if}\quad i\not=j,j+1\, , \notag \\
 (q^{}z-q^{-1}w)F_{i}(z)F_{i}(w)= F_{i}(w)F_{i}(z)(q^{-1}z-q^{}w)\,
 ,  \label{gln-com1}\\
(q^{-1}z-qw)F_{i}(z)F_{i+1}(w)= F_{i+1}(w)F_{i}(z)(z-w)\, , \notag \\
k_i^\pm(z)F_i(w)\left(k_i^\pm(z)\right)^{-1}=
\frac{q^{-1}z-qw}{z-w}F_i(w)\, , \notag \\
k_{i+1}^\pm(z)F_i(w)\left(k_{i+1}^\pm(z)\right)^{-1}=
\frac{q^{}z-q^{-1}w}{z-w}F_i(w)\, , \notag \\
k_i^\pm(z)F_j(w)\left(k_i^\pm(z)\right)^{-1}=F_j(w),
 \qquad {\rm if}\quad i\not=j,j+1\, , \notag \\
[E_{i}(z),F_{j}(w)]= \delta_{{i},{j}}\ \delta(z/w)\
(q-q^{-1})\left( k^+_{i}(z)/k^+_{i+1}(z)-k^-_{i}(w)/k^-_{i+1}(w)\right)\, \notag
\end{gather}
and the Serre relations for the  currents $E_{i}(z)$ and $F_{i}(z)$
\begin{equation}
\begin{split}
{\rm Sym}_{z_1,z_{2}}
\Big(E_{i}(z_1)E_{i}(z_2)E_{i\pm 1}(w)
&-(q+q^{-1})E_{i}(z_1)E_{i\pm 1}(w)E_{i}(z_2)+\\
&+E_{i\pm 1}(w)E_{i}(z_1)E_{i}(z_2)\Big)=0\, ,\\
\label{serre}
{\rm Sym}_{z_1,z_{2}}
\Big(F_{i}(z_1)F_{i}(z_2)F_{i\pm 1}(w)
&-(q+q^{-1})F_{i}(z_1)F_{i\pm 1}(w)F_{i}(z_2)+\\
&+F_{i\pm 1}(w)F_{i}(z_1)F_{i}(z_2)\Big)=0\, .
\end{split}
\end{equation}
The generating series $F_i(z)$, $E_i(z)$ and $k_j^\pm(z)$ are called
{\em total and Cartan currents} respectively.
Formulae \eqref{gln-com1} and \eqref{serre} should be considered as
formal series identities describing an infinite set of relations between modes of the
currents. The symbol $\delta(z)$ entering these relations is the formal series $\sum_{n\in\ZZ} z^n$.

The algebra $\Uqglnd$ in its current realization can be also obtained in the framework of
the quantum double construction choosing another type of the Borel subalgebra. It can be
constructed as the quantum double from the subalgebra $U_F$ generated by the modes of the
currents
$F_i[n]$, $k^+_j[m]$, $i=1,\ldots,N-1$, $j=1,\ldots,N$, $n\in\ZZ$
and $m\geq0$. One may easily see from the commutation relations \r{gln-com1}
that this is an subalgebra in $\Uqglnd$, but for the quantum double construction one has
to choose  for this subalgebra a coalgebraic structure different from \r{coL+} \cite{D88}:
\begin{equation}\label{coD}
\begin{split}
\Delta^{(D)}\sk{F_i(z)}&=1\ot F_i(z) +  F_i(z)\ot k^+_{i}(z)\sk{k^+_{i+1}(z)}^{-1},\\
\Delta^{(D)}\sk{k^+_i(z)}&=k^+_i(z)\ot k^+_{i}(z)\,.
\end{split}
\end{equation}
We call $U_F$ {\em a current Borel subalgebra} of $\Uqglnd$.

The dual current Borel subalgebra
$U_E\subset \Uqgln$ is generated by modes of the currents
$E_i[n]$, $k^-_j[-m]$, $i=1,\ldots,N-1$, $j=1,\ldots,N$, $n\in\ZZ$ and
$m\geq0$ with coalgebraic structure
\begin{equation*}
\begin{split}
\Delta^{(D)}\sk{E_i(z)}&= E_i(z)\ot 1 +  k^-_{i}(z)\sk{k^-_{i+1}(z)}^{-1}\ot E_i(z),\\
\Delta^{(D)}\sk{k^-_i(z)}&=k^-_i(z)\ot k^-_{i}(z).
\end{split}
\end{equation*}

In the $\LL$-operator formulation of $\Uqglnd$ we do not use the restriction \r{rest}.
In the current realization of the same algebra we shall not assume the relations
\begin{equation}\label{rest1}
k^+_i[0]k^-_i[0]=1\,,\qquad i=1,\ldots,N\,.
\end{equation}
The latter relations are standard  in the realization of the
quantum affine algebra $\Uqgln$ through a
quantum double construction \cite{D88}.
This realization
of the quantum affine algebra $\Uqgln$ can be obtained from $\Uqglnd$ by imposing the relation
\r{rest1}.

\subsection{Cartan-Weyl generators of $\Uqglnd$}

The current realization of $\Uqglnd$ uses currents corresponding to the simple roots of
$\mathfrak{gl}_N$ only. The modes of these currents can be identified with a part of the Cartan-Weyl
generators of this quantum affine algebra. Instead of the rest Cartan-Weyl generators
we will use the generating series introduced for the first time in \cite{DKh} where they were called
{\it the composed currents}.

Denote by  $\overline U_F$ an extension of the current Borel subalgebra $U_F$ formed
by the linear combinations of series, given as infinite sums of monomials
$a_{i_1}[n_1]\cdots a_{i_k}[n_k]$ with $n_1\leq\cdots\leq n_k$, and $n_1+...+n_k$
fixed, where  $a_{i_l}[n_l]$ is either $F_{i_l}[n_l]$ or $k^+_{i_l}[n_l]$.
Analogously, denote by $\overline U_E$ an extension of the dual current Borel subalgebra $U_E$ formed
by the linear combinations of series, given as infinite sums of monomials
$b_{i_1}[n_1]\cdots b_{i_k}[n_k]$ with $n_1\leq\cdots\leq n_k$, and $n_1+...+n_k$
fixed, where  $b_{i_l}[n_l]$ is either $E_{i_l}[n_l]$ or $k^-_{i_l}[n_l]$.

First we define the composed currents $F_{j,i}(t)$, $i<j$ as the series with coefficients in
$\overline U_F$. The definition of the composed currents may be written in the analytical form
\begin{equation}\label{rec-f1}
\ff_{j,i}(t)\, =\,
-\,\mathop{\rm res}\limits_{w=t}\ff_{j,a}(t)\ff_{a,i}(w)\,\frac{dw}{w}\, =\,
\mathop{\rm res}\limits_{w=t}\ff_{j,a}(w)\ff_{a,i}(t)\,\frac{dw}{w}
\end{equation}
for any $a=i+1,\ldots, j-1$. It is equivalent to the relation
\begin{equation}\label{rec-f111}
\begin{split}
\ff_{j,i}(t)&=
\oint \ff_{j,a}(t) \ff_{a,i}(w)\ \frac{dw}{w}-
\oint \frac{q^{-1}-qt/w}{1-t/w}
\;\ff_{a,i}(w) \ff_{j,a}(t)\ \frac{dw}{w}\,,\\
\ff_{j,i}(t)&=
\oint \ff_{j,a}(w) \ff_{a,i}(t)\ \frac{dw}{w}-
\oint  \frac{q^{-1}-qw/t}{1-w/t}
\;\ff_{a,i}(t) \ff_{j,a}(w)\ \frac{dw}{w}\,.\\
\end{split}
\end{equation}
In \r{rec-f111} we set $\oint \frac{dw}{w} g(w)=g_0$ for any formal series $g(w)=\sum_{n\in\ZZ}g_n z^{-n}$.

Using the  relations \r{gln-com1} on $F_i(t)$ we can calculate
the residues in \r{rec-f1} and  obtain the following expressions  for  $\ff_{j,i}(t)$, $i<j$:
\begin{equation}\label{res-in}
\ff_{j,i}(t)=(q-q^{-1})^{j-i-1}\ff_i(t) \ff_{i+1}(t)\cdots \ff_{j-1}(t)\,.
\end{equation}
For example, $F_{i+1,i}(t)=F_i(t)$, and
$F_{i+2,i}(t)=(q-q^{-1})F_i(t)F_{i+1}(t)$.
Formulas \r{res-in} prove the consistency of the defining relations for the composed currents
\r{rec-f1} or \r{rec-f111}, since they  yield to the same answers for all possible
values $i<a<j$.

Calculating formal integrals in \r{rec-f111} we obtain
the following presentation for the composed currents:
\begin{equation}\label{rec-f112a}
F_{j,i}(t)=
F_{j,a}(t)F_{a,i}[0]-q^{-1}F_{a,i}[0]F_{j,a}(t)+
(q-q^{-1})\sum_{k< 0}F_{a,i}[k]\,F_{j,a}(t)\,t^{-k}\,,
\end{equation}
\begin{equation}\label{rec-f112b}
F_{j,i}(t)=
F_{j,a}[0]F_{a,i}(t)-qF_{a,i}(t)F_{j,a}[0]+
(q-q^{-1})\sum_{k\geq 0}F_{a,i}(t)\,F_{j,a}[k]\,t^{-k}\,.
\end{equation}

Composed currents $E_{i,j}(t)$ are defined analogously as the series with coefficients which
belong to the completed current Borel subalgebra $\overline U_E$.

The fact that current generators for the quantum affine algebras form the part of the  Cartan-Weyl basis in these
algebras was proved in \cite{KT}. There exists a natural ordering in the Cartan-Weyl basis.
If the generator $e_\gamma$ corresponds to a positive root $\gamma=\alpha+\beta$, where $\alpha$
and $\beta$ are the roots, then these generators are ordered either in a way
$e_\alpha\prec e_\gamma\prec e_\beta$ or in the way
$e_\beta\prec e_\gamma\prec e_\alpha$. An important property of the Cartan-Weyl basis of a Borel
subalgebra is that the $q$-commutator of any two generators from this subalgebra, say $e_\alpha$ and $e_\beta$,
is a linear combination of monomials which contain only the products of generator $e_{\gamma_i}$
which are `between' $e_\alpha$ and $e_\beta$: $e_\alpha\prec e_{\gamma_i}\prec e_\beta$
or $e_\alpha\succ e_{\gamma_i}\succ e_\beta$.
An important application of this property is that using it one can order
arbitrary monomials of the generators.

The ordering on the Borel subalgebra can be extended to the ordering of the whole set of
Cartan-Weyl
generators corresponding to the positive and negative roots
such that the same ordering property is valid. This ordering is called `circular' or
`convex' and it
allows to order arbitrary monomials in the whole algebra \cite{EKhP}.
For the goals of our paper we will need
the following specialization of this circular ordering of the current generators in the algebra
$\Uqglnd$.

\begin{figure}[hbt]
\begin{center}
\begin{picture}(240,160) \thicklines
\put(40,135){$\cU_{q}(\fb^-)$}
\put(30,120){$\overbrace{\qquad\qquad}$}
\put(120,135){$\cU_{q}(\fb^+)$}
\put(110,120){$\overbrace{\qquad\qquad}$}
\put(0,60){\line(1,0){45}}
\put(68,60){\line(1,0){50}}
\put(140,60){\vector(1,0){50}}
\put(35,80){\line(-1,0){15}}
\put(0,80){$\cU_{E}^-$}
\put(55,78){\oval(40,60)}
\put(50,90){$\cU_{e}^-$}
\put(50,55){$\cU_{k}^-$}
\put(50,20){$\cU_{f}^-$}
\put(120,90){$\cU_{e}^+$}
\put(95,95){\oval(100,25)}
\put(145,95){\line(1,0){15}}
\put(165,90){$\cU_{e}$}
\put(120,20){$\cU_{f}^+$}
\put(95,25){\oval(100,25)}
\put(45,25){\line(-1,0){15}}
\put(15,20){$\cU_{f}$}
\put(121,55){$\cU_{k}^+$}
\put(130,41){\oval(40,57)}
\put(150,35){\line(1,0){15}}
\put(170,30){$\cU_{F}^+$}
\multiput(0,44)(5,0){17}{.}
\multiput(84,45)(3,3){11}{.}
\multiput(118,77)(5,0){24}{.}
\multiput(100,-15)(0,5){30}{.}
\put(210,95){$\displaystyle\left.\rule{0mm}{4ex}\right\}$}
\put(225,95){$\cU_E$}
\put(210,40){$\displaystyle\left.\rule{0mm}{6ex}\right\}$}
\put(225,35){$\cU_F$}
\end{picture}
\end{center}
\caption{\footnotesize{Subalgebras of $\Uqglnd$.
Vertical dotted line separates standard Borel subalgebras.
Horizontal dotted line separates current Borel subalgebras.
Horizontal solid axis shows increasing of the modes
of the current generators. Ovals signify different subalgebras in the standard and
 current Borel
subalgebras of $\Uqglnd$.}}\label{fig:subalg}
\end{figure}
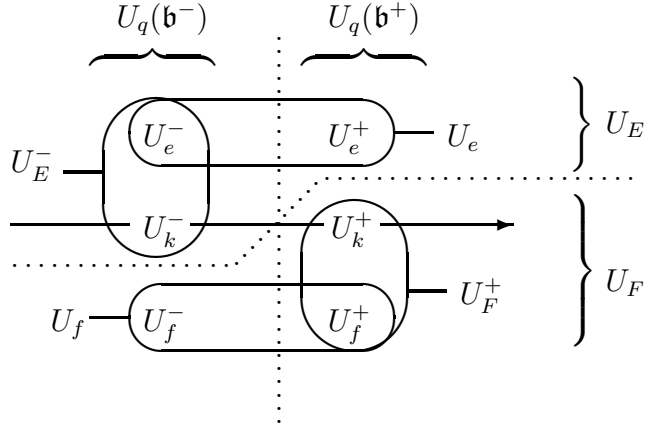

Define the intersections of the different type Borel subalgebras in $\Uqglnd$:
\begin{equation}\label{int-sec}
U^-_E=U_E \cap \Uqbm\,,\quad U^-_f= \Uqbm \cap U_F \,,\quad
U^+_F=U_F \cap \Uqbp\,,\quad U^+_e= \Uqbp \cap U_E\,.
\end{equation}
Let $U_f$ and $U_e$ be subalgebras of $\Uqglnd$ formed by the modes of the
currents $F_i(z)$ and $E_i(z)$ respectively.
Subalgebras $U^-_E$ and $U^+_F$ can be also decomposed into
subalgebras $U^-_e$, $U^-_k$ and $U^+_f$, $U^+_k$, where
\begin{equation}\label{int-sec1}
U^-_e=U_e \cap U^-_E\,,\quad U^+_f= U_f \cap U^+_F
\end{equation}
and $U^\pm_k$ are defined by the isomorphisms
$$
U_f^+\ot U^+_k \to U_F^+=U^+_f\cdot U^+_k\,,\quad
U_e^-\ot U^-_k \to U_E^-=U^-_e\cdot U^-_k\,.
$$
Different relations between these subalgebras are represented
schematically in Figure \ref{fig:subalg}.

We fix a `circular' ordering `$\prec$' on the generators of $\Uqglnd$ (see \cite{EKhP}), such that:
\begin{equation}\label{circular}
\cdots\ \prec\ U^-_k\ \prec\ U^-_f\ \prec\ U^+_f\ \prec\ U^+_k\ \prec\ U^+_e\ \prec\ U^-_e\ \prec\ U^-_k\
\prec\ \cdots \,.
\end{equation}
\begin{definition}\label{def:NO}
We will call an element $W\in\Uqglnd$ {\em normal ordered}  and
denote it as $:W:$ if it is presented as the
linear combinations of the products  the
elements $W_1\cdot W_2\cdot W_3\cdot W_4\cdot W_5\cdot W_6$ such that
\begin{equation}\label{nor-ord}
 W_1\in U^-_f\,,\quad W_2\in U^+_f\,,\quad
W_3\in U^+_k\,,\quad W_4\in U^+_e\,,\quad W_5\in U^-_e \,,\quad    W_6\in U^-_k\,.
\end{equation}
\end{definition}
We may consider standard Borel subalgebras as ordered with respect to
the circular ordering \r{circular}:
$$
\Uqbm=U^-_e\cdot U^-_k\cdot U^-_f\,,\quad \Uqbp=U^+_f\cdot U^+_k\cdot U^+_e\,.
$$
Analogous statement is valid for the current Borel subalgebras:
$$
U_F=U^-_f\cdot U^+_f\cdot U^+_k\,,\quad U_E=U^+_e\cdot U^-_e\cdot U^-_k\,.
$$

\subsection{Projections and universal off-shell Bethe vectors}

It was proved in \cite{KPT} that the subalgebras $U_f^-$
and $U_F^+$ are coideals with respect to the Drinfeld coproduct
\r{coD}
\begin{equation*}
\Delta^{(D)}(U_F^+)\subset \Uqglnd\ot U_F^+\,,\qquad
\Delta^{(D)}(U_f^-)\subset U_f^-\ot \Uqglnd\,,
\end{equation*}
and that the multiplication $m$ in $\Uqglnd$ induces an isomorphism of vector spaces
$$m: U_f^-\ot U_F^+\to U_F\,.$$
According to the general theory presented in \cite{EKhP} we  define
 projection operators $\Pfp:U_F\subset \Uqglnd \to U_F^+$ and
$\Pfm:U_F\subset \Uqglnd \to U_f^-$ by the prescriptions
\begin{equation}\label{pgln}
\begin{split}
\Pfp(\Ff_-\ \Ff_+)&=\coun(\Ff_-)\ \Ff_+, \qquad
\Pfm(\Ff_-\ \Ff_+)=\Ff_-\ \coun(\Ff_+),
\\& \text{for any}\qquad \Ff_-\in U_f^-,
\quad \Ff_+\in U_F^+ .
\end{split}
\end{equation}
It was also proved in \cite{EKhP} that
\begin{itemize}
\item[(1)] projections \r{pgln} can be extended to the
 algebra
$\overline U_F$;
\item[(2)] for any $\Ff\in \overline U_F$ with $\Delta^{(D)}(\Ff)=\sum_i \Ff'_i\otimes \Ff''_i$ we have
\begin{equation}\label{pr-prop}
\Ff=\sum_i\Pfm(\Ff''_i)\cdot \Pfp(\Ff'_i)\,.
\end{equation}
\end{itemize}
Analogously, we may define projections of the dual current Borel subalgebras $U_E$:
$\Pem : U_E \to U_E^-$ and $\Pep : U_E \to U_e^+$.

Let $\bar n=\{n_{1},n_{2},\ldots,n_{N-2},n_{N-1}\}$
be a set of non-negative integers.  Denote by  $\bar t_{\segg{\bar n}}$
 the set of formal variables:
\begin{equation}\label{set111}
\bar{t}_{\segg{\bar{n}}} = \left\{ t^{1}_{1},\ldots,t^{1}_{n_{1}};
t^{2}_{1},\ldots, t^{2}_{n_{2}};\ \ldots\ldots\ ;
t^{N-2}_{1},\ldots, t^{N-2}_{n_{N-2}};
t_1^{N-1},\ldots,t^{N-1}_{n_{N-1}}\right\}\,.
\end{equation}
The variable $t_{k}^a$ is of type $a$. If $n_a=0$ for some $a$, then
the variables of type $a$ are absent in the set \r{set111}.
Denote by $\mathcal{W}_{N}(\bar{t}_{\segg{\bar{n}}})$
the universal weight function associated with the set of  variables
\r{set111}:\footnote{In contrast to the paper \cite{FKPR-GS}, here we will not
normalize the universal weight function by the product of factorials.}
\begin{equation}\label{uwf}
\mathcal{W}_{N}(\bar{t}_{\segg{\bar{n}}})=\Pfp\left(
F_{N-1}(t_{n_{N-1}}^{N-1})\cdots F_{N-1}(t_{1}^{N-1})
  \quad \cdots\quad
F_{1}(t_{n_1}^1)\cdots F_{1}(t_{1}^1)
\right).
\end{equation}

The universal weight function  \r{uwf}
is a formal series over the ratios
 $t^b_k/t^c_l$ with $b<c$ and $t^a_i/t^a_j$ with $i<j$
 taking values in the completion ${\overline{U}}_F$.

It was proved in \cite{EKhP,KPT} that this projection satisfies the comultiplication
properties of the off-shell Bethe vectors of the hierarchical Bethe ansatz \cite{KR83}.
In \cite{KP}, a method of calculation of the projection \r{uwf}  based on the
ordering property \r{pr-prop} was proposed. Then in \cite{KhP-GLN,OPS}, this method
was applied to the quantum affine algebra $\Uqgln$. It was proved there
that the
calculation of the projection \r{uwf} produces the same hierarchical relations
for the universal off-shell Bethe vectors as were found in \cite{VT2} using combinatorial
methods in the framework of the standard hierarchical Bethe ansatz.

We call a vector $v$ \emph{a weight singular vector} if it is
annihilated by any positive mode  $E_i[n]$, $i=1,\ldots,N-1$,
$n> 0$ and is an eigenvector for  $k^+_i(t)$, $i=1,\ldots,N$
\begin{equation}\label{hwv}
\EE^{+}_{i,i+1}(t)\cdot v=0\ ,\qquad k^+_i(t)\cdot v=\lambda_i(t)\,v\,,
\end{equation}
where $\lambda_i(t)$ is a meromorphic function, decomposed as a
power series in $t^{-1}$. The $L$-operator \r{L-op}, acting on a
weight singular vector $v$, becomes upper-triangular
\begin{equation}\label{L-op-tr}
\LL^{+}_{ij}(t)\ v=0\,,\quad i>j\,,\quad \LL^{+}_{ii}(t)\ v=\lambda_i(t)\ v\,,\quad
i=1,\ldots,N\,.
\end{equation}

\begin{definition}\label{def:J}
We note $J$ the left ideal of $\Uqbp$, generated by all elements of the form
$\Uqbp\cdot E_i[n]$, $i=1,\ldots,N-1$, $n>0$ (equivalently, by all modes of $\Uqbp\cdot
{\EE}^+_{i,j}(t)$, $1\leq i<j\leq N$).\\
Equalities in $\Uqbp$ modulo element from the ideal $J$ we denote by the
symbol `$\,\simJ$'.
\end{definition}
It is clear that $W\cdot v=0$ for any element $W\in J$ and arbitrary
weight singular vector $v$.

We call {\em a (universal) transfer matrix} the trace of $\LL$-operator
\begin{equation}\label{utm}
\mathcal{T}_{N}(t)=\sum_{i=1}^N \LL_{i,i}^+(t)=
\sum_{i=1}^N\sk{ k_i^+(t)+\sum_{j=i+1}^N \FF^+_{j,i}(t)\, k^+_j(t)\, \EE^+_{i,j}(t)}.
\end{equation}
The Gauss coordinates $\FF^+_{j,i}(t)$, $\EE^+_{i,j}(t)$ coincide with the projections of the corresponding
composed currents and can be expressed through modes of the currents from subalgebras
$U_f^+$ and $U_e^+$ (see section~\ref{sect3.1}).
Note that the presentation \r{utm} of the  transfer matrix
$\mathcal{T}_{N}(t)$ is normal ordered according to the circular ordering
\r{circular} and $\mathcal{T}_{N}(t)\ \simJ\
\sum_{i=1}^N k_i^+(t)$.

The main statement of this paper is
\begin{theorem}\label{main-th}
A formal series identity is valid in $\Uqbp$
\begin{equation}\label{main-id}
\mathcal{T}_{N}(t)\cdot \mathcal{W}_{N}(\bar{t}_{\segg{\bar{n}}})\ \,
-\mathcal{W}_{N}(\bar{t}_{\segg{\bar{n}}})\cdot\tau_{N}(t;\bar{t}_{\segg{\bar{n}}})
\ \simJ\ 0
\end{equation}
if the set $\{t^i_j\}$ of the Bethe parameters satisfies the set of the universal Bethe equations,
$i=1,\ldots,N-1$, $j=1,\ldots,n_i$:
\begin{equation}\label{univ-BE}
\frac{k^+_i(t^i_j)}{k^+_{i+1}(t^i_j)}=
\prod_{m\neq j}^{n_i}\frac{q-q^{-1}t^{i}_m/t^i_j}{q^{-1}-qt^{i}_m/t^i_j}\
\prod_{m=1}^{n_{i-1}}\frac{1-t^{i-1}_m/t^i_j}{q-q^{-1}t^{i-1}_m/t^i_j}\
\prod_{m=1}^{n_{i+1}}\frac{q^{-1}-qt^{i+1}_m/t^i_j}{1-t^{i+1}_m/t^i_j}\,.
\end{equation}
We have introduced
\begin{equation}\label{univ-ev}
\tau_{N}(t;\bar{t}_{\segg{\bar{n}}})
=\sum_{i=1}^N k^+_i(t)\prod_{j=1}^{n_{i-1}}\frac{q-q^{-1}t^{i-1}_j/t}{1-t^{i-1}_j/t}
\prod_{j=1}^{n_{i}}\frac{q^{-1}-qt^{i}_j/t}{1-t^{i}_j/t}
\end{equation}
 an eigenvalue of the universal transfer matrix.

\end{theorem}

Note that the universal transfer-matrix \r{utm} is a power series on $1/t$. This is a reason for
presentation of the universal eigenvalue in the form \r{univ-ev}. By  similar reasons
the right hand side of the universal Bethe equations \r{univ-BE} is a series depending
on the ratios $t^a_m/t^i_j$.

Proofing the statement of Theorem~\ref{main-th} we will try to present the product
$\mathcal{T}_{N}(t)\cdot \mathcal{W}_{N}(\bar{t}_{\segg{\bar{n}}})$ in the
normal ordered form according to the ordering given in Definition~\ref{def:NO}.
After performing this ordering we will observe that subtraction of the ordered
product $\mathcal{W}_{N}(\bar{t}_{\segg{\bar{n}}})\cdot\tau_{N}(t;\bar{t}_{\segg{\bar{n}}})$
results only in the terms which belong to the ideal $J$.

For any weight singular vector $v$, let ${w}_V^{N}(\bar
t_{\segg{\bar n}})=\mathcal{W}_{N}(\bar t_{\segg{\bar n}})\ v$ be
the  weight function taking value in the $\Uqglnd$-module $V$
generated by $v$ and
\begin{equation}\label{b-v-sym}
{\bf w}_V^{N}(\bar t_{\segg{\bar n}})=\qsym(\bar t_{\segg{\bar
n}})\prod_{a=2}^{N}\prod_{\ell=1}^{n_{a-1}}\lambda_a(t^{a-1}_\ell) \
{w}_{V}^{N}(\bar t_{\segg{\bar n}})\end{equation}  be the corresponding modified
weight function \cite{KPT}. Here
\begin{equation*}
\qsym(\bar t_{\segg{\bar n}})=\prod_{a=1}^{N-1}\prod_{1\leq \ell<\ell'\leq n_a}
\frac{q-q^{-1}t^a_\ell/t^a_{\ell'}}{1-t^a_\ell/t^a_{\ell'}}\ .
\end{equation*}
According to \cite{KPT,KhP-GLN,VT2}
we call the modified weight function \r{b-v-sym} {\em  universal off-shell Bethe vector}.
The Theorem~\ref{main-th} has obvious
\begin{corollary}
Universal off-shell Bethe vector is an eigenvector of the universal transfer matrix
\begin{equation*}
\mathcal{T}_{N}(t)\cdot \mathbf{w}_V^{N}(\bar{t}_{\segg{\bar{n}}})
=\mathbf{w}_V^{N}(\bar{t}_{\segg{\bar{n}}})\cdot
\sk{\sum_{i=1}^n \lambda_i(t)\prod_{j=1}^{n_{i-1}}\frac{qt-q^{-1}t^{i-1}_j}{t-t^{i-1}_j}
\prod_{j=1}^{n_{i}}\frac{q^{-1}t-qt^{i}_j}{t-t^{i}_j}}
\end{equation*}
if Bethe equations of $\Uqgln$ hierarchical Bethe ansatz are satisfied:
\begin{equation*}
\frac{\lambda_i(t^i_j)}{\lambda_{i+1}(t^i_j)}=
\prod_{m\neq j}^{n_i}\frac{qt^i_j-q^{-1}t^{i}_m}{q^{-1}t^i_j-qt^{i}_m}\
\prod_{m=1}^{n_{i-1}}\frac{t^i_j-t^{i-1}_m}{qt^i_j-q^{-1}t^{i-1}_m}\
\prod_{m=1}^{n_{i+1}}\frac{q^{-1}t^i_j-qt^{i+1}_m}{t^i_j-t^{i+1}_m}\,.
\end{equation*}
\end{corollary}

\section{Proofs}\label{sect3}

We will prove the Theorem~\ref{main-th} by induction over $N$. First, we check that
the statement of the Theorem is valid in the simplest case $N=2$. Then, assuming
the correctness of the statement for the algebra $\Uqglnd$ we will prove it for
the algebra $\Uqd{N+1}$. The $\LL$-operator of the
smaller algebra $\Uqglnd$ will be embedded into right down corner of the
$\LL$-operator for the bigger algebra $\Uqd{N+1}$. This embedding is in accordance
with the Gauss decomposition \r{L-op}. The main technical tool will be a special presentation
of the universal weight function based on the main property of the projections
\r{pgln}. This allows to reduce the calculation of the ordering to the commutations
of the Gauss coordinates and total currents, which is much simpler than the exchange relations
of these coordinates with the projection \r{uwf}.

\subsection{Relation between Gauss coordinates and the currents}\label{sect3.1}
\def\f{F}
\def\e{E}

In order to perform the proof of the main Theorem we need to establish a precise relation
between {\em all}\/ Gauss coordinates and the currents. This was done partially in \cite{KPT} and
here we repeat this calculations for the sake of completeness. Note that the paper \cite{DF}
yields these relations only for currents corresponding to simple roots.

Set $S_A(B)=BA-q^{-1}AB$.
Projections of composed currents can be defined using $q$-commutators
with zero modes of the currents $F_i(t)$, $i=1,\ldots,N-1$.
We will call the operators $S_{F_i[0]}\equiv S_i$ {\it the screening operators}.
\begin{proposition}\label{com-cur2}
We have
\begin{equation}\label{pro-in}
\Pfp\sk{F_{j,i}(t)}\,=\,
S_{i}\bigl(\Pfp(F_{j,i+1}(t))\bigr)\,,\quad
\Pfm\sk{F_{j,i}(t)}\,=\,-q\,
S_{j-1}\bigl(\Pfm(F_{j-1,i}(t))\bigr)\,,
\quad i<j-1\,.
\end{equation}
\end{proposition}
\noindent
{\it Proof.}\enspace
We apply the projection $\Pfp$ to both sides of the relation
\r{rec-f112a} with $a=i+1$
and  the projection $\Pfm$ to both sides of the relation
\r{rec-f112b} with $a=j-1$.
The modes $F_i[k]$ with $k<0$ belong to $U^-_f$ and the modes
$F_{j-1}[k]$ with $k\geq0$ belong to $U^+_f$. Hence, due to
formulae~\r{pgln}, the projections $\Pfpm$ both kill  the semi-infinite sum in
the right hand side of \r{rec-f112a} and \r{rec-f112b}, and we get
\begin{equation}\label{pro1}
\begin{split}
\Pfp\sk{F_{j,i}(t)}\,&{}=\,
\Pfp\sk{F_{j,i+1}(t)F_i[0]-q^{-1}\,F_i[0]F_{j,i+1}(t)}\,={}\\[4pt]
&{}=\,\Pfp\sk{S_{i}\sk{F_{j,i+1}(t)}}\,=\,
S_{i}\sk{\Pfp\sk{F_{j,i+1}(t)}}
\end{split}
\end{equation}
and
\begin{equation}\label{pro2}
\begin{split}
\Pfm\sk{F_{i,j}(t)}\,&{}=\,-q\,
\Pfm\sk{\,F_{j-1,i}(t)F_{j,j-1}[0]-q^{-1}F_{j,j-1}[0]F_{j-1,i}(t)}\,={}\\[4pt]
&{}=\,-q\,\Pfm\sk{S_{j-1}\sk{F_{j-1,i}(t)}}\,=\,-q\,
S_{j-1}\sk{\Pfm\sk{F_{j-1,i}(t)}}
\end{split}
\end{equation}
To get the last equalities we use the fact proved in \cite{KP} that
the projection $\Pfpm$ commutes with the screening operators $S_{i}$:
$\Pfpm\bigl(S_{i}(W)\bigr)=S_{i}\bigl(\Pfpm(W)\bigr)$
for any $W\in U_f$.
\hfill$\square$

The screening operators also relate the Gauss coordinates of the $\LL$-operators.
\begin{lemma}\label{5.6}
We have
\begin{equation}\label{cor1}
(q-q^{-1})\FF^\pm_{j,i}(t)\,=\,S_{i}\sk{\FF^\pm_{j,i+1}(t)}\,,\qquad i<j-1\,.
\end{equation}
\end{lemma}
\noindent
{\it Proof}\enspace of this Lemma was given in \cite{KPT}\footnote{See also the proof of the
analogous  Lemma~\ref{5.6e} below.} and is based on the
commutation relation between matrix elements $\LL^+_{i,i+1}(s)$ and
$\LL^\pm_{i+1,j}(t)$.
One should consider the coefficients at the zero power of the spectral parameter $s$ in these
relations
and take into account that $\LL^+_{i,i+1}[0]=F_{i}[0]k^+_{i+1}[0]$ and
$\LL^+_{i+1,i+1}[0]=k^+_{i+1}[0]$.
\qed

 Proposition~\ref{pro-in} and Lemma~\ref{5.6} imply the following
\begin{proposition}\label{idenGC}
We have
\begin{equation}\label{ide-pr+}
\Pfp\sk{F_{j,i}(t)}\,=\,(q-q^{-1})^{j-i-1}\FF^+_{j,i}(t)\,,\qquad i<j-1\,,
\end{equation}
\begin{equation}\label{ide-pr-}
\Pfm\sk{F_{j,i}(t)}\,=\,-(q-q^{-1})^{j-i-1}\sk{\FF^-_{j,i}(t)+\sum_{\ell=1}^{j-i-1}
(-1)^\ell\!\!\!\! \sum_{j>i_{\ell}>\cdots>i_1>i}\FF^-_{i_1,i}(t)\cdots \FF^-_{j,i_\ell}(t)
}\,.
\end{equation}
\end{proposition}
\noindent{\it Proof.}\enspace
First equality \r{ide-pr+} was proved in \cite{KPT} using
  induction with respect to $j-i$ from the formula
$\Pfp\bigl(F_{i+1,i}(t)\bigr)=\FF^+_{i+1,i}(t)$ \cite{DF}.
Here we shall prove \r{ide-pr-}. We apply the projection $\Pfm$ to both sides of the
relation \r{rec-f112a} to obtain
\begin{equation}\label{pmgc1}
\Pfm\sk{F_{j,i}(t)}=S_i\sk{\Pfm\sk{F_{j,i+1}(t)}}+(q-q^{-1})\ \Pfm\sk{F_{i+1,i}(t)}\cdot
\Pfm\sk{F_{j,i+1}(t)}\,.
\end{equation}
Using this relation recursively
 and the Lemma~\ref{5.6}, we get
\begin{equation}\label{pmgc2}
\Pfm\sk{F_{j,i}(t)}+(q-q^{-1})^{j-i-1}\FF^-_{j,i}(t)+\sum_{\ell=i+1}^{j-1}(q-q^{-1})^{j-\ell}\
\FF^-_{\ell,i}(t)\cdot \Pfm\sk{\FF_{j,\ell}(t)}=0\,.
\end{equation}
>From the identity
\begin{equation}
\Pfm\sk{F_{i+1,i}(t)}= -\FF^-_{i+1,i}(t)
\end{equation}
one proves that equality \r{ide-pr-} is a solution of this
recurrence relation, and coincide with it for $i=j-1$.
\qed

Proceeding in analogous way, we may relate the projections of the dual composed currents
$\Pepm(E_{i,j}(t))$ with Gauss coordinates $\EE^\pm_{i,j}(t)$, but here we shall need only
the relations between different dual Gauss coordinates or analog of the Lemma~\ref{5.6} for
$\EE^\pm_{i,j}(t)$. Set $\tS_A(B)=AB-qBA$ and denote $\tS_{E_i[0]}\equiv \tS_i$.
\begin{lemma}\label{5.6e}
We have
\begin{equation}\label{cor1e}
(q-q^{-1})\EE^\pm_{i,j}(t)\,=\,\tS_{i}\sk{\EE^\pm_{i+1,j}(t)}\,,\qquad i<j-1\,.
\end{equation}
\end{lemma}
\noindent
{\it Proof.}\enspace
Let us consider the commutation relations between the following matrix elements of $\LL$-operators
\begin{equation*}
\begin{split}
(t-s)[\LL^\pm_{j,i+1}(t),\LL^-_{i+1,i}(s)]&=(q-q^{-1})\sk{
t\LL^-_{i+1,i+1}(s)\LL^\pm_{j,i}(t) - s\LL^\pm_{i+1,i+1}(t)\LL^-_{j,i}(s)}\,,\\
(qt-q^{-1}s)\LL^-_{i+1,i+1}(s)\LL^\pm_{j,i+1}(t)&= (t-s) \LL^\pm_{j,i+1}(t)\LL^-_{i+1,i+1}(s)
+(q-q^{-1})s \LL^\pm_{i+1,i+1}(t)\LL^-_{j,i+1}(s)\,.
\end{split}
\end{equation*}
Choosing the coefficients at the zero power of the spectral parameter $s$ in these relations
and taking into account that $\LL^-_{i+1,i}[0]=-k^-_{i+1}[0]E_i[0]$ and $\LL^-_{i+1,i+1}[0]=k^-_{i+1}[0]$
we obtain
\begin{equation}\label{imp1}
(q-q^{-1})\LL^\pm_{j,i}(t)=
E_i[0]\LL^\pm_{j,i+1}(t)-q\LL^\pm_{j,i+1}(t)E_i[0]=\tS_i\sk{\LL^\pm_{j,i+1}(t)}\,.
\end{equation}
In order to obtain \r{cor1e} we shall use an explicit expression for the matrix elements of
the $\LL$-operator \r{GE1} in terms of the Gauss coordinates. The relation \r{imp1} implies
\r{cor1e} for $j=N$ due to the commutativity of $k_N^\pm(t)$ and $E_i[0]$ for $i=1,\ldots,N-2$.
Next, the relations \r{imp1} for $j=N-1$ and \r{cor1e} for $j=N$ imply \r{cor1e} for $j=N-1$
due to the commutativity of the Gauss coordinates $k^\pm_{N-1}(t)$,  $k^\pm_{N}(t)$ and
$\FF^\pm_{N,N-1}(t)$ with $E_i[0]$ for $i=1,\ldots,N-3$. The statement of the Lemma follows by
induction over $j$.
\qed

\subsection{Basic notations}

Let $\bar\lll$ and $\bar\rr$ be two collections of nonnegative integers
satisfying a set of inequalities
\begin{equation}\label{set23}
\lll_a\leq\rr_a\,,\quad a=1,\ldots,N-1\,.
\end{equation}
Denote by $\meg{\bar\lll}{\bar\rr}$  a set of segments which contain  positive
integers $\{\lll_a+1,\lll_a+2,\ldots,\rr_a-1,\rr_a\}$
including $\rr_a$ and excluding $\lll_a$. The length of each segment is equal
to $\rr_a-\lll_a$.

For a given set $\meg{\bar\lll}{\bar\rr}$ of segments  we denote by
$\bar t_{\meg{\bar\lll}{\bar\rr}}$  the sets of variables
\begin{equation}\label{t-coll}
\bar t_{\seg{\bar\rr}{\bar\lll}} =
\{t^{1}_{\lll_{1}+1},\ldots,t^{1}_{\rr_{1}};
t^{2}_{\lll_{2}+1},\ldots,t^{2}_{\rr_{2}};\ldots;
t^{N-1}_{\lll_{N-1}+1},\ldots,t^{N-1}_{\rr_{N-1}} \}.
\end{equation}
{}For any $a=1,\ldots,N-1\,$ we denote the sets of variables corresponding to the  segments
$\meg{\lll_a}{\rr_a}=\{\lll_a+1,\lll_a+2,\ldots,\rr_a\}$
as $\bar{t}^a_{\meg{{\lll_a}}{\rr_a}}\
=\{t^{a}_{\lll_{a}+1},\ldots,t^{a}_{\rr_{a}} \}$.
All the variables in $\bar{t}^a_{\meg{{\lll_a}}{\rr_a}}$ have type $a$. For the segments
$\meg{\lll_a}{\rr_a}=\meg{0}{n_a}$ we use the shorten  notations
$\bar{t}_{\meg{\bar 0}{\bar n}}\equiv\bar{t}_{\megg{\bar n}}$ and
$\bar{t}^a_{\meg{0}{n_a}}\equiv\bar{t}^a_{\megg{n_a}}$.

For a collection of variables $\bar t_{\seg{\bar\rr}{\bar\lll}}$ we
consider the ordered products of the currents
\begin{equation}\label{tFFFl}
\F(\bar t_{\seg{\bar\rr}{\bar\lll}})=\!\!\prod_{N-1\geq a\geq
1}^{\longleftarrow} \sk{\prod_{\rr_a\geq\ell>
\lll_a}^{\longleftarrow} F_{a}(t^a_{\ell})}=
F_{N-1}(t^{N-1}_{\rr_{N-1}})\cdots
F_{1}(t^{1}_{\rr_{1}})\cdots F_{1}(t^1_{\lll_1+1})\,,
\end{equation}
where the series $F_a(t)\equiv F_{a+1,a}(t)\,$ are
 defined by  \r{DF2}.
 As particular cases, we have
$\,\F(\bar t^a_{\seg{\rr_a}{\lll_a}})=
F_{a}(t^a_{r_a})\cdots F_{a}(t^a_{l_a+2})
F_{a}(t^a_{l_a+1})$.

The product \r{tFFFl}
is a formal series over the ratios
 $t^b_k/t^c_l$ with $b<c$ and $t^a_i/t^a_j$ with $i<j$
 taking values in the algebra  ${{U}}_F$.

 Symbol $\mathop{\prod}\limits^{\longleftarrow}_a A_a$ (resp.
$\mathop{\prod}\limits^{\longrightarrow}_a A_a$)
 means the  ordered products of
noncommutative entries $A_a$, such that $A_a$ is on the right (resp. on the left)
 from $A_b$ for $b>a$:
\begin{equation*}
\mathop{\prod}\limits^{\longleftarrow}_{j\geq a\geq i} A_a = A_j\,A_{j-1}\,
\cdots\, A_{i+1}\,A_i\,,\qquad
\mathop{\prod}\limits^{\longrightarrow}_{i\leq a\leq j} A_a = A_i\,A_{i+1}\,
\cdots\, A_{j-1}\,A_j\,.
\end{equation*}

Consider the permutation group $S_n$ and its action on the formal series of $n$ variables
defined for the elementary transpositions $\sigma_{i,i+1}$ as follows
\begin{align*}
\pi(\sigma_{i,i+1})G(t_1, \dots,t_i, t_{i+1},\dots, t_n) =
\frac{q^{-1}-q\,t_i/t_{i+1}}{q-q^{-1}\,t_i/t_{i+1}}\ G(t_1, \dots,t_{i+1}, t_i, \dots, t_n),
\end{align*}
The $q$-depending factor in this formula is chosen in such a way that each product
$F_a(t_n)\cdots$ $F_a(t_1)$
is invariant under this action. Summing the action over all the group of permutations
we obtain the operator $\tSym_u=\frac{1}{n!}\sum\limits_{\sigma\in S_n}\pi(\sigma)$ acting as
follows\footnote{Normalization of the $q$-symmetrization used here differs from the one used
in the papers \cite{KhP-GLN,OPS} by the combinatorial factor $\frac{1}{n!}$.}
\begin{align}\label{qss}
\tSym_{\ \bar t}\  G(\bar t) =\frac{1}{n!}
\sum\limits_{\si \in S_n}\prod\limits_{\substack{\ell<\ell'\\ \si(\ell)>\si(\ell')}}
\frac{q^{-1}-q\,t_{\si(\ell')}/t_{\si(\ell)}}
{q-q^{-1}\,t_{\si(\ell')}/t_{\si(\ell)}}\  G(^\sigma t).
\end{align}
The product is taken over all pairs $(\ell, \ell')$, such that conditions $\ell
< \ell'$ and $\si(\ell) > \si(\ell')$ are satisfied simultaneously.

We call the operator  $\tSym_u$  a {\em $q$-symmetrization}.
The operator $\tSym_u$ is the group average
with respect to the action $\pi$, so that
\begin{equation}\label{sym*}
\tSym_{\ \bar t}\ \tSym_{\ \bar t}\ (\cdot)=\ \tSym_{\ \bar t}\ (\cdot)\,.
\end{equation}

An important property of $q$-symmetrization is the relation
\begin{align}
 \tSym_{(t_1,\ldots,t_n)}
  =\frac{s!(n-s)!}{n!}\sum_{\sigma\in
    S^{(s)}_n}\pi(\sigma)\;\tSym_{(t_1,\ldots,t_s)}\tSym_{(t_{s+1},\ldots,t_n)}\,,
    \label{sym_div}
\end{align}
where $s\in\seg{n}{0}$ is fixed and the sum is taken over the subset
\begin{align*}
 S^{(s)}_n=\{\sigma\in S_n \mid \sigma(1)<\ldots<\sigma(s)\,;\ \sigma(s+1)<\ldots<\sigma(n)\}\,.
\end{align*}

Denote by $S_{[\bar\lll,\bar\rr]} =
S_{[\lll_{1},\rr_{1}]}\times \cdots \times S_{[\lll_{N-1},\rr_{N-1}]}$
 the direct product of the groups $S_{[\lll_{a},\rr_{a}]}$ permuting
integer numbers $\lll_{a}+1,\ldots, \rr_{a}$.
The $q$-sym\-me\-tri\-za\-tion over the whole set of variables
$\bar t_{\seg{\bar\rr}{\bar\lll}}$ is defined by the formula
\begin{equation}\label{qsr}
\tSym_{\ \bar t_{\seg{\bar\rr}{\bar\lll}}} \ G(\bar
t_{\seg{\bar\rr}{\bar\lll}})= \sum_{\si\in
S_{[\bar\lll,\bar\rr]}}\prod_{1\leq a\leq N-1}
\left(\frac{1}{(\rr_a-\lll_a)!}
\prod_{\substack{\ell<\ell'\\ \si^a(\ell)>\si^a(\ell')}}
\frac{q^{-1}-q^{}\,t^a_{\si^a(\ell')}/t^a_{\si^a(\ell)}}
{q^{}-q^{-1}\,t^a_{\si^a(\ell')}/t^a_{\si^a(\ell)}}\right)\ G(^\si \bar
t_{\seg{\bar\rr}{\bar\lll}})\,,
\end{equation}
where  the set $^\si \bar t_{\seg{\bar\rr}{\bar\lll}}$
is defined  as
\begin{equation}\label{sigmat}
^\si \bar t_{\seg{\bar\rr}{\bar\lll}} =
\{t^{1}_{\si^1(\lll_{1}+1)},\ldots,t^{1}_{\si^1(\rr_{1})};
t^{2}_{\si^2(\lll_{2}+1)},\ldots,t^{2}_{\si^2(\rr_{2})};\ldots;
t^{N-1}_{\si^{N-1}(\lll_{N-1}+1)},\ldots,t^{N-1}_{\si^{N-1}(\rr_{N-1})} \}.
\end{equation}

We say that the series $G(\bar t_{\seg{\bar\rr}{\bar\lll}})$ is $q$-symmetric, if
it is invariant under the action $\pi$ of each group
$S_{[\lll_{a},\rr_{a}]}$ with respect to the permutations of the variables
$t^a_{\lll_a+1},\ldots,t_{\rr_a}$ for $a=1,\ldots,N-1$:
\begin{equation}\label{exa3}
\tSym_{\ \bar t_{\seg{\bar\rr}{\bar\lll}}}
G(\bar t_{\seg{\bar\rr}{\bar\lll}})=
G(\bar t_{\seg{\bar\rr}{\bar\lll}})\,.
\end{equation}
The $q$-symmetrization $G(\bar t_{\seg{\bar\rr}{\bar\lll}})=
\tSym_{\bar t_{\seg{\bar\rr}{\bar\lll}}} Q(\bar t_{\seg{\bar\rr}{\bar\lll}})$
of any series $Q(\bar t_{\seg{\bar\rr}{\bar\lll}})$ is a $q$-symmetric series,
which follows from \r{sym*}.

Let $\bar s=\{s_1,\ldots,s_{N-1}\}$ be a set of nonnegative integers satisfying
$l_a\leq s_a\leq r_a$ for $a=1,\ldots,N-1$. The set of integers $\bar s$ divides the set
of the variables $\bar t_{\seg{\bar r}{\bar l}}$ into two subsets $\bar t_{\seg{\bar r}{\bar s}}\cup
\bar t_{\seg{\bar s}{\bar l}}$.

Using the property of the projections \r{pr-prop}
we can present any product of the currents in a normal ordered form
(in the sense of definition \ref{def:NO}):
\begin{equation}\label{dec-ff1m}
\begin{array}{c}
\ds \F(\bar t_{\seg{\bar\rr}{\bar\lll}})=
\sum_{\lll_{N-1}\leq \ss_{N-1}\leq \rr_{N-1}}\cdots \sum_{\lll_1\leq \ss_1\leq \rr_1}
\ \  \prod_{1\leq a\leq N-1}  \frac{(\rr_a-\lll_a)!}{(\ss_a-\lll_a)!(\rr_a-\ss_a)!}\times \\ [10mm]
\ds \times\   \tSym_{\ \bar t_{\seg{\bar\rr}{\bar\lll}}}
\left(Z_{\bar\ss}({\bar t}_{\seg{\bar\rr}{\bar\lll}})
\ds\  \Pfm\sk{\F(\bar t_{\seg{\bar\rr}{\bar\ss}})}\cdot
\Pfp\sk{\F(\bar t_{\seg{\bar\ss}{\bar\lll}})}\right),
\end{array}
\end{equation}
where
\begin{equation}\label{Zserm}
Z_{\bar\ss}({\bar t}_{\seg{\bar\rr}{\bar\lll}})=\prod_{a=1}^{N-2}\ \
\prod_{\substack{ \ss_a < \ell\leq \rr_a \\
  \lll_{a+1}< \ell' \leq \ss_{a+1}}} \frac{q-q^{-1}\
t^{a}_{\ell}\,/\,t^{a+1}_{\ell'}}{1-t^{a}_{\ell}\,/\,t^{a+1}_{\ell'}}\,.
\end{equation}
Equality \r{dec-ff1m} was proved in \cite{KhP-GLN} and the proof is based on the
current coproduct property \r{coD} and the exchange  relations between currents.

\subsection{Special presentation of the universal weight function}

Let ${\bar \mm}=\{\mm_1,\ldots,\mm_{N-1}\}$
be a collection of the non-negative integers
satisfying the admissibility condition
\begin{equation}\label{m-non-inc}
\mm_1\geq\mm_{2}\geq\mm_{3}\ldots\geq \mm_{N-1}\geq\mm_{N}=0\,.
\end{equation}
 We define a series depending on the set of the variables
$\bar t_{\segg{\bar m}}$
 of the form
\begin{equation}\label{t-rat-X}
\tilde X(\bar t_{\segg{\bar\mm}})=
\prod_{a=1}^{N-2}
V(t^{a+1}_{\mm_{a+1}},\ldots,t^{a+1}_{1};
t^{a}_{\mm_{a}},\ldots,t^{a}_{\mm_a-\mm_{a+1}+1} )\,.
\end{equation}
where the rational series $V(\cdot;\cdot)$ is given by the formulae
\begin{equation}\label{rat-Y}
\begin{split}
\ds \tilde V(t^2_k,\ldots,t^2_1;t^1_k,\ldots,t^1_1)&=\ds \prod_{m=1}^k
\left(\frac{1}{1-t^1_m/t^2_m}
\prod_{m'=m+1}^{k}\frac{q-q^{-1}t^1_{m'}/t^2_m}{1-t^1_{m'}/t^2_m}\right)\\
&=\ds \prod_{m=1}^k\left(\frac{1}{1-t^1_m/t^2_m}
\prod_{m'=1}^{m-1}\frac{q-q^{-1}t^1_m/t^2_{m'}}{1-t^1_m/t^2_{m'}}\right)\ .
\end{split}
\end{equation}

Define a normalized ordered  product of the composed currents:
\begin{equation}\label{dstring}
\tSgs_{\bar\mm}(\bar t_{\segg{\bar\mm}})=\tilde X(\bar t_{\segg{\bar\mm}})
\prod^{\longleftarrow}_{N\geq a > 1}\sk{\frac{1}{ (\mm_{a-1}-\mm_{a})! }
\prod^{\longleftarrow}_{\mm_1-\mm_{a}\geq \ell>\mm_1-\mm_{a-1}}
F_{a,1}(t^1_\ell)}\,.
\end{equation}
This ordered product was called {\it the dual string} in the paper \cite{FKPR-GS}. Denote
the negative projections of the $q$-symmetrized dual strings \r{dstring} as follows
\begin{equation}\label{coef-m}
\Ags_{m_1, m_2,\ldots,m_{N-1}}(\bar t_{\segg{\bar\mm}})
=\Pfm\sk{\tSym_{\ \bar t_{\segg{\bar\mm}}}\
\sk{\tSgs_{\bar\mm}(\bar t_{\segg{\bar\mm}})}}\,.
\end{equation}
Denote by $\Bgs_{n_1,\ldots,n_{N-1}}(\bar t_{\segg{\bar n}})$ the elements of $U^-_f$
defined by the recursive relations
\begin{equation}\label{rec1}
\tSym_{\ \bar t_{\segg{\bar n}}}\!\!
\sk{
\mathop{\sum_{n_1\geq m_1\geq 0}\cdots \sum_{n_{N-1}\geq m_{N-1}\geq 0}}
\limits_{m_1\geq \cdots \geq m_{N-1}}
 Z_{\bar m}(\bar t_{\segg{\bar n}})
\Bgs_{n_1-m_1,\ldots,n_{N-1}-m_{N-1}}(\bar t_{\seg{\bar n}{\bar m}})
\cdot \Ags_{m_1,\ldots,m_{N-1}}(\bar t_{\segg{\bar m}})}\!\!=0.
\end{equation}
It was proved in \cite{FKPR-GS} that the coefficients $\Bgs_{n_1,\ldots,n_{N-1}}(\bar t_{\segg{\bar n}})$
defined by \r{rec1} are non-zero only iff $n_1\geq n_2\geq\cdots \geq n_{N-1}\geq 0$ and
can be defined uniquely by means of \r{rec1} from the initial condition
$\Bgs_{0,\ldots,0}(\bar t_{\segg{\bar n}})=1$. In particular,
\begin{equation}\label{partB}
\Bgs_{\underbrace{\scriptstyle 1,\ldots,1}_{m\ {\rm times}}\!\!,0,\ldots,0}(t^1,\ldots,t^{N-1})=
-\Ags_{\underbrace{\scriptstyle 1,\ldots,1}_{m\ {\rm times}}\!\!,0,\ldots,0}(t^1,\ldots,t^{N-1})=
-\prod_{a=1}^{m-1}\frac{1}{1-t^a/t^{a+1}}\ \Pfm\sk{F_{m+1,1}(t^1)}.
\end{equation}

Denote by
$\bar t_{\segg{\bar\ss'}}$ the following collection of the formal variables
\begin{equation*}
\bar t_{\segg{\bar\ss'}}=\{t^2_{1},\ldots,t^2_{\ss_2};
\ldots;t^{N-1}_{1},\ldots,t^{N-1}_{\ss_{N-1}}\}\,.
\end{equation*}
excluding the variables of type 1.
We formulate without proof the following
\begin{proposition}\label{spelpr} \cite{FKPR-GS}
There is a formal series identity
\begin{equation}\label{spelpr1}
\begin{split}
\mathcal{W}_{N}(\bar t_{\segg{n}})&=\sum_{\bar\ss} \prod_{a=1}^{N-1}\frac{n_a!}{s_a!}\ \
  \tSym_{\ \bar t_{\segg{n}}} \Big( Z_{\bar\ss}(\bar t_{\segg{\bar n}})\ \times \\
&\times
\Bgs_{n_1-s_1,\ldots,n_{N-1}-s_{N-1}}(\bar t_{\seg{\bar n}{\bar\ss}})
\cdot
 \mathcal{W}_{N-1}(\bar t_{\segg{\ss'}}) \cdot F_1(t^1_{\ss_1})\cdots F_1(t^1_{1})\Big)\,.
\end{split}
\end{equation}
\end{proposition}

In this paper we will need the following Corollary of this Proposition.
Let $\bar t_{\segg{\bar n'}_m}$ be the following collection of formal variables
\begin{equation}\label{set-m-min}
\bar t_{\segg{\bar n'}_m}=\{t^2_1,\ldots,t^2_{n_2-1};\ldots;
t^m_1,\ldots,t^m_{n_m-1};t^{m+1}_1,\ldots,t^{m+1}_{n_{m+1}};\ldots; t^{N-1}_1,\ldots,t^{N-1}_{n_{N-1}}\}\,.
\end{equation}
Note that $\bar t_{\segg{\bar n'}_1}\equiv \bar t_{\segg{\bar n'}}$.
\begin{corollary}\label{sp-el-pr}
\begin{equation}\label{sef1}
\begin{split}
\Pfp\sk{\F(\bar t_{\segg{\bar n}})}=&\Pfp\sk{\F(\bar t_{\segg{\bar n'}_1})}
\cdot \F(\bar t^1_{\segg{n_1}})-\sum_{m=1}^{N-1} \prod_{a=1}^m (n_a)
\ \tSym_{\ \bar t_{\segg{\bar n}}}
(\Pfm(F_{m+1,1}(t^1_{n_1}))\times\\
&\quad\times\  \Pfp(\F(\bar t_{\segg{\bar n'}_m}))
\cdot  \F(\bar t^1_{\segg{n_1-1}})\cdot \Zfun_m(\bar t_{\segg{\bar n}})  ) + \mathbb{W}\,,
\end{split}
\end{equation}
where
\begin{equation}\label{Zfun}
\Zfun_m(\bar t_{\segg{\bar n}})=   \prod_{a=1}^{m-1}\left(
\frac{1}{1-t^{a}_{n_{a}}/t^{a+1}_{n_{a+1}}}
\prod_{j=1}^{n_{a+1}-1}
\frac{{q-q^{-1}t^{a}_{n_{a}}/t^{a+1}_j}}{1-t^{a}_{n_{a}}/t^{a+1}_j}
\right)
\prod_{j=1}^{n_{m+1}}\frac{{q-q^{-1}t^{m}_{n_{m}}/t^{m+1}_j}}{1-t^{m}_{n_{m}}/t^{m+1}_j}
\end{equation}
and the terms $\mathbb{W}$ in \r{sef1} are such that $\Pfp\sk{:\mathcal{T}_{N}(t)\cdot \mathbb{W}:}=0$.
\end{corollary}
Recall that $\F(\bar t^1_{\segg{n_1}})=F_1(t^1_{n_1})\cdots F_1(t^1_{1})$ and
$\F(\bar t^1_{\segg{n_1-1}})=F_1(t^1_{n_1-1})\cdots F_1(t^1_{1})$.

The first term in the right hand side of \r{sef1} corresponds to the term with
all $s_a=n_a$ in \r{spelpr1}. Each of the terms in the summation over $m$ in \r{sef1}
corresponds
to the following values of $s_m$ in the general formula  \r{spelpr1}:
$s_1=n_1-1,\ldots,s_m=n_m-1$ and $s_{m+1}=n_{m+1},\ldots,s_{N-1}=n_{N-1}$. The corresponding
elements $\Bgs_{1,\ldots,1,0,\ldots,0}$  are given by \r{partB}, which brings in \r{sef1}
the product of the rational factors $(1-t^{a}_{n_{a}}/t^{a+1}_{n_{a+1}})^{-1}$. Other rational
factors are given by the series $Z_{\bar\ss}(\bar t_{\segg{\bar n}})$ for these particular
values of $\bar s$.

The general structure of the terms $\mathbb{W}$ which are not presented explicitly in the right hand side
of \r{sef1} can be described as follows.
The structure of the coefficients
$\Ags_{m_1, m_2,\ldots,m_{N-1}}(\bar t_{\segg{\bar\mm}})$ \r{coef-m} implies that these
terms will have on the left the negative projections of the string
$\tSgs_{\bar\mm}(\bar t_{\segg{\bar\mm}})$ which contains, at least,  the product of two
currents $F_{c_1,1}$ and $F_{c_2,1}$ or the product of the several negative projections of the
strings of type \r{coef-m}. Since the projection of the string can be always factorized to the product of the
projection of the currents \cite{KhP-GLN}, the general structure of the terms $\mathbb{W}$ will be
$\mathbb{W}=\sum\Pfm\sk{F_{c_1,1}}\cdot
\Pfm\sk{F_{c_2,1}} \cdot  \mathbb{W}'$.
The elements $\mathbb{W}'$ are some elements of $U_F$ and their
exact structure is unimportant. The reason why $\Pfp\sk{:\mathcal{T}_{N}(t)\cdot \mathbb{W}:}=0$
will be explained in the next subsection.

Note that the identity \r{spelpr1} can be proved directly using only the ordering relations
\r{dec-ff1m} and the rules of calculations of the negative projections from the
product of currents.
Indeed, for an arbitrary product of currents $\F(\bar t_{\segg{\bar
n}})$, these
ordering relations
can be written in the form
\begin{equation}\label{rrrs}
\Pfp\sk{\F(\bar t_{\segg{\bar n}})}=\F(\bar t_{\segg{\bar n}})-\sum
\Pfm\sk{\F'}\cdot \Pfp\sk{\F''}\,,
\end{equation}
where the number of currents in the product $\F''$ is less than in the original
product $\F(\bar t_{\segg{\bar n}})$. Thus, one can replace recursively the positive projection
$\Pfp\sk{\F''}$ by the right hand side of the relation \r{rrrs} up to
 the obvious identity $\Pfp\sk{F_i(t)}=F_i(t)-\Pfm\sk{F_i(t)}$ valid for arbitrary
simple current $F_i(t)$. Calculating the negative projections $\Pfm\sk{\F'}$
to obtain the projections of the strings of type \r{coef-m}, we can prove the
identity  \r{spelpr1} by brute force calculations. The technique of the generating
series developed in \cite{FKPR-GS} yields more elegant way of proving this and many other similar
identities.

\begin{example}
Let us present an example of the general formula \r{spelpr1} in the case $N=3$, $n_1=2$ and $n_2=2$.
To reduce the formula we will use shorthand notations
$\Pfpm(\cdot)=[\cdot]^\pm$. We also denote
$t^2_i=s_i$ and $t^1_i=t_i$ and $\tSym$ below will be the $q$-symmetrization over variables $t_i$ and
$s_i$.
\begin{equation*}
\begin{split}
&\big[F_2(s_2)F_2(s_1)F_1(t_2)F_1(t_1)\big]^+=\big[F_2(s_2)F_2(s_1)\big]^+F_1(t_2)F_1(t_1)\\
&\quad-2\,\tSym\sk{[F_1(t_2)]^- \big[F_2(s_2)F_2(s_1)\big]^+ F_1(t_1)
\prod_{j=1}^2 \frac{q^{-1}t_2-qs_j}{t_2-s_j}  }\\
&\quad-4\,\tSym\sk{[F_{3,1}(t_2)]^- [F_2(s_1)]^+\, F_1(t_1)\, \frac{s_2}{s_2-t_2}
\frac{q^{-1}t_2-qs_1}{t_2-s_1}  }\\
&\quad+\left\{4\,\tSym\sk{ \sk{ [F_1(t_2)]^-\,[F_1(t_1)]^- -\frac{1}{2}\big[F_1(t_2)F_1(t_1)\big]^-}
\big[F_2(s_2)F_2(s_1)\big]^+
\prod_{i,j=1}^2 \frac{q^{-1}t_i-qs_j}{t_i-s_j}  }\right.\\
&\quad+4\,\tSym\left( \left(
\frac{s_2}{s_2-t_2}\ [F_{3,1}(t_2)]^-\,[F_1(t_1)]^-+
\frac{s_2}{s_2-t_1}\frac{q^{-1}t_2-qs_2}{t_2-s_2}\ [F_{1}(t_2)]^-\,[F_{3,1}(t_1)]^-\right.\right.\\
&\qquad\qquad \left.\left.-\frac{s_2}{s_2-t_2}\ \big[F_{3,1}(t_2)F_1(t_1)\big]^-
\right) [F_2(s_1)]^+
\prod_{i=1}^2 \frac{q^{-1}t_i-qs_1}{t_i-s_1}  \right)\\
&\left.\quad+4\,\tSym\sk{ \sk{ [F_{3,1}(t_2)]^-[F_{3,1}(t_1)]^- -\frac{1}{2}\big[F_{3,1}(t_2)F_{3,1}(t_1)\big]^-}
\frac{s_1}{s_1-t_1}\frac{s_2}{s_2-t_2} \frac{q^{-1}t_2-qs_1}{t_2-s_1}
}\right\}\,.
\end{split}
\end{equation*}
\end{example}

The terms in curly brackets correspond to the term $\mathbb{W}$ in \r{sef1}.

\subsection{The action of $\LL_{a,b}^+(t)$ onto $\Pfm\sk{F_{c,d}(t)}$}
\label{actL-F}

\begin{definition}\label{def:I}
Let $I$ be the right ideal of $\Uqglnd$, generated by all elements of the form
$F_i[n]\cdot \Uqbp$ such that $i=1,\ldots,N-1$ and  $n<0$. We denote equalities modulo
elements from the ideal $I$ by the symbol `\,$\simI$'.
\end{definition}

\begin{proposition}\label{aLFm}
We have an equivalence
\begin{equation}\label{aLF1}
\LL_{a,b}^+(t)\cdot \FF^-_{c,d}(s)\ \simI\  \delta_{a,c}\ \frac{(q-q^{-1})s}{s-t}\ \LL^+_{d,b}(t)\,.
\end{equation}
\end{proposition}

One of our technical tools will be the rule of commuting the negative projections of the
composed currents $\Pfm\sk{F_{c,d}(t')}$ with matrix elements of the fundamental $\LL$-operator
$\LL_{a,b}^+(t)$. We need a result of this calculation only modulo elements from the ideal $I$
and call this as \emph{action} of $\Pfm\sk{F_{c,d}(t')}$ onto $\LL_{a,b}^+(t)$.
Due to the relation \r{ide-pr-},
in order to calculate the action of the matrix elements
$\LL_{a,b}^+(t)$ onto $\Pfm\sk{F_{c,d}(t')}$ one has to calculate first the action of
the matrix elements  $\LL_{a,b}^+(t)$ onto Gauss coordinates $\FF^-_{c,d}(t')$.

{\it Proof}\/ of Proposition~\ref{aLFm} will be done considering  each fixed $a$.

Fix $a$ and consider $c<a$. This case is simple.
Formulas \r{ide-pr-} can be inverted
to express the Gauss coordinates $\FF^-_{c,d}(s)$ in terms of the modes of the currents
$F_d[n_d],\ldots, F_{c-1}[n_{c-1}]$. But the $\LL$-operator modes
$\LL^+_{a,b}[n]$ simply commute with these current modes and so
$\LL_{a,b}^+(t)\cdot \FF^-_{c,d}(s)=\FF^-_{c,d}(s)\cdot \LL_{a,b}^+(t)\in I$.

The case when $a=N$ is also simple. For this choice, we have also $b=N$  and
$\LL^+_{N,N}(t)\equiv k^+_N(t)$ commutes with the Gauss coordinates
$\FF^-_{c,d}(s)$ for $c=2,\ldots,N-1$.
It means that  $\LL_{N,N}^+(t)\cdot \FF^-_{c,d}(s)\sim 0$ for $c<N$.
Let $c=N$. Taking into account that $\FF^-_{N,d}(s)=\LL^-_{d,N}(s)\LL^-_{N,N}(s)^{-1}$ and
the commutation relation
\begin{equation}\label{sim2}
\LL^+_{N,N}(t)\LL^-_{d,N}(s)=\frac{qt-q^{-1}s}{t-s}\LL^-_{d,N}(s)\LL^+_{N,N}(t)-
\frac{(q-q^{-1})s}{t-s}\LL^+_{d,N}(t)\LL^-_{N,N}(s)
\end{equation}
we prove the statement of the Proposition for $a=N$.

Let $a<N$ and $c>a$. First consider the case when $b>c>a$. We have
\begin{equation*}
\LL_{a,b}^+(t)\cdot \FF^-_{c,d}(s)=\FF^+_{b,a}(t)\FF^-_{c,d}(s)k^+_b(t)+
\sum_{j=b+1}\FF^+_{j,a}(t)\FF^-_{c,d}(s)k^+_j(t)E^+_{b,j}(t)
\end{equation*}
due to the commutativity of the Gauss coordinates $k^+_b(t)$ and
$k^+_j(t)E^+_{b,j}(t)$ with modes of the currents $F_d[n_d],\ldots, F_{c-1}[n_{c-1}]$
or with Gauss coordinates $\FF^-_{c,d}(s)$. The statement of the proposition follows
from the lemma.
\begin{lemma}\label{ser-ac}
For $b>c>d>a$ and $b>c>a\geq d$
\begin{equation*}
\FF^+_{b,a}(t)\FF^-_{c,d}(s)\simI F_{b,a}(t)\FF^-_{c,d}(s)\in I\,.
\end{equation*}
\end{lemma}
\noindent
{\it Proof}\/ is based on the commutation relations of the composed currents
$F_{b,a}(t)$ and $F_{c,d}(s)$.
They are
\begin{equation*}
\begin{split}
\ds \frac{t-s}{qt-q^{-1}s}\ F_{b,a}(t)\,F_{c,d}(s)&=\ds\frac{q^{-1}t-qs}{t-s}\ F_{c,d}(s)\,F_{b,a}(t),\qquad
d<a\,,\\
\ds F_{b,a}(t)\,F_{c,d}(s)&=\ds\frac{q^{-1}t-qs}{t-s}\ F_{c,d}(s)\,F_{b,a}(t),\qquad d=a\,,\\
\ds F_{b,a}(t)\,F_{c,d}(s)&= F_{c,d}(s)\,F_{b,a}(t),\ \ \qquad\qquad\qquad d>a\,,
\end{split}
\end{equation*}
and they take into account the Serre relations \r{serre} (see details in
Appendix A of the paper \cite{KhP-GLN}).
The  product $F_{b,a}(t)\,F_{c,d}(s)$ has
no poles for $d\geq a$ and has first order pole at the point $t=s$  in the case $d<a$, but the residue
at this point of this product is zero. It means that commuting negative projections
of the current $\Pfm\sk{F_{c,d}(s)}$ through the total currents $F_{b,a}(t)$ no higher currents
will be created and the result of commutation will belong to the right ideal $I$. Because
of the relation between negative Gauss coordinates and the negative projections
of the composed currents given by \r{ide-pr-} the same statement
will be true for the Gauss coordinates. \qed

Next we consider the cases when $c>a$ and $c\geq b$.
The statement of Proposition~\ref{aLFm}
will be proved by a double induction over $c$ starting from $c=N$ and over $b$ starting from $c$.
Let $c=N$ and $b=N$. Then
 using again the fact
$\FF^-_{N,d}(s)=\LL^-_{d,N}(s)k^-_{N}(s)^{-1}$ and the commutation relations
\begin{equation}\label{sim1}
\LL^+_{a,N}(t)\LL^-_{d,N}(s)=\frac{t-s}{q^{-1}t-qs}\LL^-_{d,N}(s)\LL^+_{a,N}(t)+
\frac{(q-q^{-1})s}{q^{-1}t-qs}\LL^-_{a,N}(t)\LL^+_{d,N}(t)
\end{equation}
we obtain the inclusion $\LL^+_{a,N}(t)\FF^-_{N,d}(s)\in I$.
Before considering other cases we prove the following lemma.
\begin{lemma}\label{muta}
For $b<c$ and arbitrary $a<b$ and $d<c$ we have
\begin{equation*}
\begin{split}
\LL^+_{a,b}(t)&\cdot \LL^-_{d,c}(s) +\frac{(q-q^{-1})^2ts}{(t-s)^2}\
\LL^+_{d,c}(t)\cdot \LL^-_{a,b}(s) \simI 0\,,\qquad a\neq d\,,\\
\LL^+_{a,b}(t)&\cdot \LL^-_{a,c}(s) \simI 0\,,\qquad a=d\,.
\end{split}
\end{equation*}
\end{lemma}
\noindent
{\it Proof.}\/ The case of $a=d$ follows from the  relation
\begin{equation*}
\LL^+_{a,b}(t)\LL^-_{a,c}(s)=\frac{t-s}{qt-q^{-1}s}\LL^-_{a,c}(s)\LL^+_{a,b}(t)+
\frac{(q-q^{-1})s}{qt-q^{-1}s}\LL^-_{a,b}(s)\LL^+_{a,c}(t)\,.
\end{equation*}
The case $a<d$ follows from two  relations
\begin{equation*}
[\LL^+_{a,b}(t),\LL^-_{d,c}(s)]=\frac{(q-q^{-1})}{t-s}\sk{s\ \LL^+_{a,c}(t)\LL^-_{d,b}(s)-
t\ \LL^-_{a,c}(s)\LL^+_{d,b}(t)}\,,
\end{equation*}
\begin{equation*}
[\LL^+_{a,c}(t),\LL^-_{d,b}(s)]=\frac{(q-q^{-1})t}{t-s}\sk{ \LL^-_{d,c}(s)\LL^+_{a,b}(s)-
 \LL^+_{d,c}(t)\LL^-_{a,b}(s)}\,.
\end{equation*}
The case $a>d$ can be proved analogously. \qed

Return to the proof of Proposition~\ref{aLFm}.
Keep $c=N$ and consider $b=N-1$. Then
\begin{equation}\label{ca1}
\begin{split}
&\LL^+_{a,N-1}(t)\FF^-_{N,d}(s)= \LL^+_{a,N-1}\LL^-_{d,N}(s) k^-_{N}(s)^{-1}\simI\\
&\quad\simI
-\frac{(q-q^{-1})ts}{(t-s)^2}\LL^+_{d,N}(t)\LL^-_{a,N-1}(s)k^-_{N}(s)^{-1}=\\
&\quad=-\frac{(q-q^{-1})ts}{(t-s)^2}
\LL^+_{d,N}(t)\Big(\FF^-_{N-1,a}(s)  k^-_{N-1}(s)+ \FF^-_{N,a}(s)  k^-_{N}(s)\EE^-_{N-1,N}(s)
\Big)k^-_{N}(s)^{-1},
\end{split}
\end{equation}
where we used Lemma~\ref{muta}.  Now the first term in the right hand side of \r{ca1} belongs to the ideal $I$
because the second index of $\LL^+_{d,N}(t)$ is bigger than the first index of $\FF^-_{N-1,a}(s)$.
The second term corresponds to the case $c=b=N$ considered above,
and thus also belongs to $I$.
Reducing $b$ and using
the Lemma~\ref{muta} we proved the statement for all $b<N=c$.

Let now $c=N-1$. For the negative Gauss coordinate $\FF^-_{N-1,d}(s)$ we can use the formula
\begin{equation}\label{f-sim}
\FF^-_{N-1,d}(s)=\sk{\LL^-_{d,N-1}(s)- \LL^-_{d,N}(s)\EE^-_{N-1,N}(s) }k^-_{N-1}(s)^{-1}.
\end{equation}
To prove that the  product $\LL^+_{a,N-1}(t)\LL^-_{d,N-1}(s)\in I$ we can use the same arguments
as for the product $\LL^+_{a,N}(t)\LL^-_{d,N}(s)$. The fact that
$\LL^+_{a,N-1}(t)\LL^-_{d,N}(s)\in I$ was already proved above. Continuing we check
that $\LL^+_{a,b}(t)\LL^-_{d,N-1}(s)\in I$ for all $b<N-1$. For general $c$ we have to
use instead of the formulae \r{f-sim} the relation
\begin{equation}\label{f-gen}
\FF^-_{c,d}(s)=\LL^-_{d,c}(s)k^-_{c}(s)^{-1}+ \LL^-_{d,c+1}(s)X_{c+1}+\cdots + \LL^-_{d,N}(s)X_{N}\,,
\end{equation}
where the explicit form of the elements $X_j\in \Uqbm$ is not important.

At last we have to check the case $a=c<N$. The case $a=c=N$ was considered above.
According to \r{f-gen} the consideration of $\LL^+_{a,b}(t) \FF^-_{a,d}(s)$ reduces to the analysis
of the product $\LL^+_{a,b}(t) \LL^-_{d,a}(s) k^-_a(s)^{-1}$.
We have
\begin{equation*}
[\LL^+_{a,b}(t),\LL^-_{d,a}(s)]=\frac{(q-q^{-1})}{t-s}\sk{t\ \LL^-_{d,b}(s)\LL^+_{a,a}(t)-
s\ \LL^+_{d,b}(t)\LL^+_{a,a}(t)}\,.
\end{equation*}
Since $\LL^-_{a,a}(s) k^-_a(s)^{-1}=1+\sum_{j=a+1}^N \FF^-_{j,a}(s)k^-_{j}(s)\EE^-_{a,j}(s)k^-_a(s)^{-1}$ and
$\LL^+_{d,b}(t)\FF^-_{j,a}(s)\in I$ the statement of the Proposition~\ref{aLFm} is proved.\qed

\begin{corollary}\label{aLFm1}
We have an equivalence
\begin{equation}\label{aLF11}
\LL_{a,b}^+(t)\cdot \Pfm\sk{F_{c,d}(s)}\simI \delta_{a,c}\ \frac{(q-q^{-1})^{c-d-1}s}{t-s}\ \LL^+_{d,b}(t)\,.
\end{equation}
\end{corollary}
\noindent
{\it Proof.}\/ Let us apply the matrix element $\LL_{a,b}^+(t)$ to both side of \r{ide-pr-}.
The first term gives the right hand side of \r{aLF11}. Other terms produce  products
of the Kronecker's symbols $\delta_{a,i_1}\,\delta_{d,i_2}\,\delta_{i_1,i_3}\cdots \delta_{i_{\ell-1},c}$
which are zero due to the restriction $d<i_1<\cdots<i_\ell<c$. If $\ell=1$ then
$\delta_{a,i_1}\,\delta_{d,c}=0$ since $d<c$. \qed

Let us explain why $\Pfp\sk{:\mathcal{T}_{N}(t)\cdot \mathbb{W}:}=0$, where $\mathbb{W}$
are the terms not shown explicitly
in the right hand side of \r{sef1}. Due to Corollary~\ref{aLFm1} the action of $\LL_{a,b}^+(t)$
onto the product of two negative projections of the currents $\Pfm\sk{F_{c_1,1}}\cdot
\Pfm\sk{F_{c_2,1}}$ is proportional to the product of delta-symbols:
$\delta_{a,c_1}\delta_{1,c_2}$. But  since $c_2>1$ this is zero modulo elements of the right ideal
$I$, which obviously satisfies $\Pfp(I)=0$.

\section{Ordering of the universal objects}\label{sect4}

The proof of  main Theorem~\ref{main-th} consists of a detailed analysis of the circular
ordering of the product of the transfer matrix and of the  universal Bethe vectors expressed
in terms of the current generators of the quantum affine algebra $\Uqglnd$.
 In the next two subsections we will perform such an analysis to the case $N=2$
and prove main Theorem~\ref{main-th} for the algebra $\Uqgldva$.
Then we go on by induction over $N$.

In what follows, besides the right ideal $I$
and the left ideal $J$ introduced in Definition~\ref{def:I} and
Definition~\ref{def:J}, we will also use the following ideal $K$.
\begin{definition}\label{def:K}
We denote by $K$ the two-sided $\Uqglnd$ ideal
  generated by the elements which have at least one
arbitrary mode $k^-_i[n]$, $i=1,\ldots,N$, $n\leq0$, of the negative Cartan
current $k^-_i(t)$.\\
Equalities in $\Uqglnd$ modulo element of the ideal $K$ are denoted by
the symbol `$\,\simK$'.\\
Equalities in $\Uqglnd$ modulo the right ideal $I$, the left ideal $J$
and the two-sided ideal $K$ will be denoted by the
symbol `$\,\simm$'.
\end{definition}

\subsection{Ordering for $\Uqgldva$}

The algebra $\Uqgldva$ is generated by the modes of the Gauss coordinates $k_1^\pm(t)$, $k_2^\pm(t)$,
$\EE_{12}^\pm(t)$, $\FF_{21}^\pm(t)$ in the $\LL$-operator realization or
by the modes of the  currents $k_1^\pm(t)$, $k_2^\pm(t)$,
$E_{1}(t)$, $F_{1}(t)$ in the current realization. The standard quantum affine algebra $\Uqgltwo$ can be obtained
from $\Uqgldva$ by imposing the restriction $k_i^+[0]\,k_i^-[0]=1$, $i=1,2$.
To simplify further formulas we shall not use index of the single simple root in the
notation of the Gauss coordinates and the currents, that is we use the following identification:
$\EE_{12}^\pm(t)\equiv \EE^\pm(t)$, $\FF_{21}^\pm(t)\equiv \FF^\pm(t)$,
$E_{1}(t)\equiv E(t)$, $F_{1}(t)\equiv F(t)$. Let $\psi^\pm(t)=k_1^+(t)k_2^+(t)^{-1}$.

The universal transfer matrix for $\Uqgldva$ is given by the relation
\begin{equation*}
\mathcal{T}_2(t)=\LL^+_{11}(t)+\LL^+_{22}(t)\ \mbox{ with }\
\LL^+_{11}=k^+_1(t)+ \FF^+(t)k^+_2(t)\EE^+(t)\,,\ \LL^+_{22}(t)=
k^+_2(t)\,,
\end{equation*}
while a universal weight function is a projection $\Pfp\sk{\F(\bar t)}=\Pfp\sk{F(t_n)\cdots F(t_1)}$. We avoid to
use upper index in the notation of the formal variables $t_j$.

\begin{proposition}\label{prop4.2} There is a formal series equality in the algebra $\Uqgldva$
\begin{equation}\label{lhs3}
\begin{split}
\mathcal{T}_2(t)\cdot \Pfp\sk{\F(\bar t)}\
&\simm\ \Pfp\sk{\F(\bar t)}\sk{\prod_{i=1}^n\frac{q^{-1}-qt_i/t}{1-t_i/t}\ k^+_1(t)+
\prod_{i=1}^n\frac{q-q^{-1}t_i/t}{1-t_i/t}\ k^+_2(t)}+\\
&+ n\ \overline{\rm Sym}_{\ \bar t}
\sk{\FF^+(t)k^+_2(t)\  F(t_{n})\cdots F(t_{2})\
\frac{(q-q^{-1})t_1/t}{1-t_1/t}\  \psi^+(t_1)  }\\
&- n\ \overline{\rm Sym}_{\ \bar t}
\sk{\FF^+(t)k^+_2(t)\ F(t_{n-1})\cdots F(t_{1})\
\frac{(q-q^{-1})t_n/t}{1-t_n/t}  }\,.
\end{split}
\end{equation}
\end{proposition}

\noindent
{\it Proof}\/ of Proposition~\ref{prop4.2} is based on the special
presentation of the universal weight function given by the Corollary~\ref{sp-el-pr}.
In this case we have
\begin{equation}\label{sef2}
\Pfp\sk{\F(\bar t)}= F(t_{n})\cdots F(t_{1})- n\ \tSym_{\ \bar t}
\sk { \Pfm\sk{F(t_{n})}\cdot F(t_{n-1})\cdots F(t_{1})} + \mathbb{W}\,,
\end{equation}
where $\mathbb{W}$ are the terms which have on the left the product of at least  two negative projections
of the currents $F(t)$. Since $\Pfm\sk{F(t)}=-\FF^-(t)$ these terms can be equally described as
having on the left the product of at least of two negative Gauss coordinates $\FF^-(t)$.
As was explained above and as we will see explicitly below, the terms $\mathbb{W}$ are
characterized by the property that $\mathcal{T}_2(t)\cdot \mathbb{F}\in I$.

We will order the product of $\mathcal{T}_2(t)$ and each summand in the right hand side of
\r{sef2} separately. For the ordering of the first term we use the relation
\begin{equation*}
[\EE^+(t),F(t_1)]= \frac{(q-q^{-1})t_1}{t-t_1}(\psi^+(t_1)-\psi^-(t_1))\,,
\end{equation*}
so that
\begin{equation}\label{lhs1}
\begin{split}
&\skk{k^+_1(t)+k^+_2(t)+\FF^+(t)k^+_2(t)\EE^+(t)}\cdot F(t_n)\cdots F(t_1)\ =\\
&\quad=\ F(t_n)\cdots F(t_1)\sk{\prod_{i=1}^n\frac{q^{-1}t-qt_i}{t-t_i}\ k^+_1(t)+
\prod_{i=1}^n\frac{qt-q^{-1}t_i}{t-t_i}\ k^+_2(t)}\\
&\qquad+ n\ \overline{\rm Sym}_{\ \bar t}
\sk{\FF^+(t)k^+_2(t)\ F(t_{n})\cdots F(t_{2})\
\frac{(q-q^{-1})t_1}{t-t_1}\  \psi^+(t_1)  } \\
&\qquad- n\ \overline{\rm Sym}_{\ \bar t}
\sk{\FF^+(t)k^+_2(t)\ F(t_{n})\cdots F(t_{2})
\frac{(q-q^{-1})t_1}{t-t_1}\  \psi^-(t_1)  } \\
&\qquad+\ \FF^+(t)k^+_2(t)\ F(t_n)\cdots F(t_{1})\  \EE^+(t)\,.
\end{split}
\end{equation}
Here we used the exchange relations of the Cartan currents $k^+_i(t)$ and the total currents
$F(t_j)$ and $q$-symmetrization in  second and  third terms of the right hand side of \r{lhs1}
appears when commuting Cartan currents $\psi^\pm(t_i)$ with total currents $F(t_j)$.

Observe that the last term in \r{lhs1} belongs to the ideal $J$ and the next to last term
 belongs to the ideal $K$. At this point
we benefit from the absence of the restriction \r{rest1} in the algebra $\Uqgldva$. Otherwise
we would have
to take into account the zero modes of the currents $\psi^-(t_1)$ if we
were considering the standard
quantum affine algebra $\Uqgltwo$.

According to \r{dec-ff1m}
\begin{equation*}
F(t_n)\cdots F(t_1)\simI \Pfp\sk{F(t_n)\cdots F(t_1)}
\end{equation*}
and equality \r{lhs1} implies the equivalence
\begin{equation}\label{lhs2}
\begin{split}
&\skk{k^+_1(t)+k^+_2(t)+\FF^+(t)k^+_2(t)\EE^+(t)}\cdot F(t_n)\cdots F(t_1)\ \simm\\
&\quad\simm\ \Pfp\sk{F(t_n)\cdots F(t_1)}\sk{\prod_{i=1}^n\frac{q^{-1}-qt_i/t}{1-t_i/t}\ k^+_1(t)+
\prod_{i=1}^n\frac{q-q^{-1}t_i/t}{1-t_i/t}\ k^+_2(t)}\\
&\qquad+ n\ \overline{\rm Sym}_{\ \bar t}
\sk{\FF^+(t)k^+_2(t)\ F(t_{n})\cdots F(t_{2})\
\frac{(q-q^{-1})t_1/t}{1-t_1/t}\  \psi^+(t_1)  }\,.
\end{split}
\end{equation}

Consider now the ordering
of the product $\mathcal{T}_2(t)\sk { \Pfm\sk{F(t_{n})}\cdot F(t_{n-1})\cdots F(t_{1})}$.
To perform this   we will use the following specialization of the Proposition~\ref{aLFm}.
For arbitrary element $\mathbb{X}\in U_F$  and $m\geq 1$ we have
\begin{equation}\label{acLF1}
\begin{split}
\LL_{11}^+(t)\cdot \Pfm\sk{F(t_m)\cdots F(t_1)}\cdot \mathbb{X}&\simI 0\,,\\
\LL_{22}^+(t)\cdot \Pfm\sk{F(t_m)\cdots F(t_1)}\cdot \mathbb{X}&\simI \delta_{m1}
\frac{(q-q^{-1})t_1/t}{1-t_1/t}
\ \FF^+(t)\,k^+_2(t)\cdot \mathbb{X}\,.
\end{split}
\end{equation}

Applying \r{acLF1} to the second term in the right hand side of \r{sef2} we obtain
\begin{equation}\label{lhs444}
\mathcal{T}_2(t)\cdot \Pfm\sk{F(t_{n})}\cdot F(t_{n-1})\cdots F(t_{1})\simI
\frac{(q-q^{-1})t_n/t}{1-t_n/t}
\ \FF^+(t)\,k^+_2(t)\cdot F(t_{n-1})\cdots F(t_{1})\,.
\end{equation}
Since any equality modulo ideal $I$ implies the related equality modulo
all the ideals $I$, $J$ and $K$,
the relations \r{lhs2} and \r{lhs444} imply the statement of Proposition~\ref{prop4.2}.\qed

\begin{corollary}\label{cor4.2}
Relation \r{lhs3} can be considered as the equality
\begin{equation*}
\begin{split}
\mathcal{T}_2(t)\cdot \Pfp\sk{\F(\bar t)}\
&\simJ\ \Pfp\sk{\F(\bar t)}\sk{\prod_{i=1}^n\frac{q^{-1}-qt_i/t}{1-t_i/t}\ k^+_1(t)+
\prod_{i=1}^n\frac{q-q^{-1}t_i/t}{1-t_i/t}\ k^+_2(t)}+\\
&+ n\ \overline{\rm Sym}_{\ \bar t}
\sk{\Pfp\sk{\FF^+(t)k^+_2(t)\  F(t_{n})\cdots F(t_{2})}\
\frac{(q-q^{-1})t_1/t}{1-t_1/t}\  \psi^+(t_1)  }\\
&- n\ \overline{\rm Sym}_{\ \bar t}
\sk{\Pfp\sk{\FF^+(t)k^+_2(t)\ F(t_{n-1})\cdots F(t_{1})}\
\frac{(q-q^{-1})t_n/t}{1-t_n/t}  }\,.
\end{split}
\end{equation*}
in $\Uqbp$  modulo elements of the ideal $J$.
\end{corollary}
\noindent
{\it Proof.}\/ Left hand side of the equality \r{lhs3} belongs to $\Uqbp$ while the right hand side
does not. We can imposing  projection $\Pfp$ onto this right hand side to cancel all the terms
which belongs to the ideal $I$.\qed

\subsection{Proof of Theorem~\ref{main-th} for $\Uqgldva$}

Let us compare the last two lines in \r{lhs3}. They contain so called
`unwanted' terms. In order to cancel them  we will use the following
properties of $q$-symmetrization. For any formal series $G(t_1,\ldots,t_n)$ on $n$ formal
variables $t_i$ we have
\begin{equation}\label{sym-last}
n\ \overline{\rm Sym}_{\ \bar t}\ G(t_1,\ldots,t_n)=
\sum_{m=1}^n \prod_{j=m+1}^n
\frac{q-q^{-1}t_m/t_j}{q^{-1}-qt_m/t_j}\
\overline{\rm Sym}_{\ \bar t\setminus t_m}\ G(t_1,\ldots,t_{m-1},t_{m+1},\ldots,t_n,t_m)
\end{equation}
and
\begin{equation}\label{sym-first}
n\ \overline{\rm Sym}_{\ \bar t}\ G(t_1,\ldots,t_n)=
\sum_{m=1}^n \prod_{j=1}^{m-1}
\frac{q-q^{-1}t_j/t_m}{q^{-1}-qt_j/t_m}\
\overline{\rm Sym}_{\ \bar t\setminus t_m}\ G(t_m,t_1,\ldots,t_{m-1},t_{m+1},\ldots,t_n)\,,
\end{equation}
where $q$-symmetrization in the right hand sides of this formal series identities
runs over $(n-1)$ variables $\bar t\setminus t_m=\{t_1,\ldots,t_{m-1},t_{m+1},\ldots,t_m\}$.
Note that formulas \r{sym-last} and \r{sym-first} are particular cases of the property
\r{sym_div} for $s=n$ and $s=1$, respectively.

Using formulas \r{sym-last} and \r{sym-first} we may write the difference of unwanted terms
in \r{lhs3} as a sum
\begin{equation*}
\begin{split}
n\sum_{m=1}^n \overline{\rm Sym}_{\ \bar t\setminus t_m}&
\left( \FF^+(t)k^+_2(t)\ F(t_n)\cdots F(t_{m+1}) F(t_{m-1}) \cdots F(t_1)\ \frac{(q-q^{-1})t_m/t}{1-t_m/t}\right.\\
&\left.\sk{\prod_{j=m+1}^n \frac{q-q^{-1}t_m/t_j}{q^{-1}-qt_m/t_j}\ \psi^+(t_m)-
\prod_{j=1}^{m-1} \frac{q-q^{-1}t_j/t_m}{q^{-1}-qt_j/t_m}   }\right).
\end{split}
\end{equation*}
Each term in this sum should be cancelled separately. This is possible if the Bethe parameters
satisfy the universal Bethe equations \cite{ACDFR}
\begin{equation}\label{un-BE1}
\psi^+(t_m)=\frac{k^+_1(t_m)}{k^+_{2}(t_m)}=
\prod_{j\neq m}^{n}\frac{q-q^{-1}t_j/t_m}{q^{-1}-qt_j/t_m}\,,\quad m=1,\ldots,n\,.
\end{equation}
Thus, Theorem~\ref{main-th} in the case $N=2$ is proved. \qed

\subsection{General case}

\subsubsection{Preliminary exchange relations}

To consider the general case of Theorem~\ref{main-th},
we need the embedding $\emb: \Uqd{N} \hookrightarrow \Uqd{N+1}$, defined
by the rule\footnote{We omit writing explicitly superscript `$+$' of the $\LL$-operator and their
Gauss coordinate, assuming that they are always from the standard positive Borel subalgebra
$\Uqbp$.}
\begin{equation}\label{embed}
\emb\sk{\LL^{[N]}_{i,j}(t)}=\LL_{i+1,j+1}(t)\,,\quad i,j=1,\ldots,N\,.
\end{equation}

Let $\mathcal{T}'_N(t)=\sum_{i=2}^{N+1}L_{ii}(t)$ be the universal transfer matrix for the
algebra $\Uqd{N}$ and $\mathcal{W}'_{N}(\bar t_{\segg{\bar n'}})$ be the universal weight
function, where $\bar n'$ is a set of the positive integers $\{n_2,\ldots,$ $n_N\}$ and
$\bar t_{[\bar n']}$ is an associated set of the formal variables:
\begin{equation*}
\bar t_{\segg{\bar n'}}=\{t^2_{1},\ldots,t^2_{n_2};\ldots; t^N_{1},\ldots,t^N_{n_N}\}\,.
\end{equation*}

Using the result of  Corrolary~\ref{sp-el-pr} we present the
$\Uqd{N+1}$-universal weight function
$\mathcal{W}_{N+1}(\bar{t}_{\segg{\bar{n}}})$ in the form
\begin{equation}\label{gc1}
\mathcal{W}_{N+1}(\bar{t}_{\segg{\bar{n}}})=
\mathcal{W}'_{N}(\bar{t}_{\segg{\bar{n}'}})\cdot \F(\bar
t_{\segg{n_1}})-\mathbb{S}\,,
\end{equation}
where $\mathbb{S}$ contains the sum of terms as in the right hand side of \r{sef1} and
the redundant terms $\mathbb{W}$ such that $:\mathcal{T}_{N+1}(t)\cdot
\mathbb{W}:\ \simI 0$.

We consider  the product
\begin{equation}\label{gc2}
\mathcal{T}_{N+1}(t)\cdot \mathcal{W}_{N+1}(\bar{t}_{\segg{\bar{n}}})\ =\
\mathcal{T}_{N+1}(t)\cdot \mathcal{W}'_{N}(\bar{t}_{\segg{\bar{n}'}})\cdot
\F(\bar t_{\segg{n_1}})-\mathcal{T}_{N+1}(t)\cdot\mathbb{S}\,.
\end{equation}
Since $\mathcal{T}_{N+1}(t)=\LL_{11}(t) + \mathcal{T}'_{N}(t)$ we rewrite the
first term in the right hand side of \r{gc2} modulo terms from the ideal $I$:
\begin{equation}\label{gc3}
\mathcal{T}'_{N}(t)\cdot \mathcal{W}'_{N}(\bar{t}_{\segg{\bar{n}'}})\cdot
\F(\bar t_{\segg{n_1}}) + \LL_{11}(t)\cdot
\F(\bar t_{\segg{\bar n}})\,.
\end{equation}
To obtain \r{gc3} we use the result of the following
\begin{lemma}\label{L11}
\begin{equation}\label{gc4}
\LL_{11}(t)\cdot \mathcal{W}'_{N}(\bar{t}_{\segg{\bar{n}'}})\cdot
\F(\bar t_{\segg{n_1}})\simI \LL_{11}(t)\cdot
\F(\bar t_{\segg{\bar n}})\,.
\end{equation}
\end{lemma}
\noindent {\it Proof.}\/
The universal Bethe vector $\mathcal{W}'_{N}(\bar{t}_{\segg{\bar{n}'}})$ for the algebra $\Uqglnd$
embedded into $\Uqd{N+1}$ by \r{embed} depends on the modes of the currents $F_2(t),\ldots,F_N(t)$ and
is given by the projection $\Pfp\sk{\F(\bar{t}_{\segg{\bar{n}'}})}$ of the product of these currents.
According to the ordering rules \r{dec-ff1m} we may replace this projection by the product
of the total currents $\F(\bar{t}_{\segg{\bar{n}'}})$ subtracting the terms which have on the left
the negative projection of the currents $\Pfm\sk{F_{c,d}(t)}$ with $2\leq d<c\leq N+1$. According to
\r{aLF11} the product $\LL_{11}(t)\Pfm\sk{F_{c,d}(t')}\simI 0$. This proves the statement of the Lemma.
\qed

Now we apply the result of  Proposition~\ref{aLFm} to the second term in the
right hand side of \r{gc2}. We have
\begin{equation}\label{gc6}
\begin{split}
\mathcal{T}_{N+1}(t)\cdot\mathbb{S}\ \simI\
&\tSym_{\ \bar t_{\segg{\bar n}}}
\left(\sum_{m=1}^{N}\ \LL_{1,m+1}(t)\cdot \Pfp\sk{\F(\bar t_{\segg{\bar n'}_m})}
\cdot \F(\bar t^1_{\segg{n_1-1}})\times \right.\\
&\left.\qquad \times \prod_{a=1}^m (n_a)\
\frac{(q-q^{-1})^m t^1_{n_1}/t}{1-t^1_{n_1}/t}\ \Zfun_m(\bar t_{\segg{\bar n}})
\right)
\end{split}
\end{equation}
where the sets of the formal variables  $\bar t_{\segg{\bar n'}_m}$ are defined
by \r{set-m-min} and a rational series $\Zfun_m(\bar t_{\segg{\bar n}})$  is
defined by \r{Zfun}. We can simplify the right hand side of \r{gc6} replacing
the projection $\Pfp\sk{\F(\bar t_{\segg{\bar n'}_m})}$  by the product of the
currents $\F(\bar t_{\segg{\bar n'}_m})$ due to  the following
\begin{lemma}\label{al2}
\begin{equation}\label{gc7}
\LL_{1,m+1}(t)\cdot \Pfp\sk{\F(\bar t_{\segg{\bar n'_m}})}\simI
\LL_{1,m+1}(t)\cdot \F(\bar t_{\segg{\bar n'_m}})\,.
\end{equation}
\end{lemma}
\noindent
{\it Proof}\/ of this Lemma is analogous to the proof of Lemma~\ref{L11}.\qed

Using the explicit form of the matrix elements $\LL_{1,1}(t)$  and $\LL_{1,m+1}(t)$ in terms of the
Gauss coordinates
\begin{equation*}
\LL_{1,1}(t)=k^+_{1}(t)+\sum_{j=2}^{N+1}
\FF^+_{j,1}(t)k^+_{j}(t)\EE^+_{1,j}(t)\,,
\end{equation*}
\begin{equation*}
\LL_{1,m+1}(t)=\FF^+_{m+1,1}(t)k^+_{m+1}(t)+\sum_{j=m+2}^{N+1}
\FF^+_{j,1}(t)k^+_{j}(t)\EE^+_{m+1,j}(t)\,,
\end{equation*}
we can present the product in the left hand side of \r{gc2} as the sum
of the terms modulo ideal $I$
\begin{equation}\label{gc8}
\begin{split}
&\mathcal{T}_{N+1}(t)\cdot \mathcal{W}_{N+1}(\bar{t}_{\segg{\bar{n}}})\ \simI\
\mathcal{T}'_{N}(t)\cdot \mathcal{W}'_{N}(\bar{t}_{\segg{\bar{n}'}})\cdot
\F(\bar t_{\segg{n_1}})+k^+_1(t)\cdot \F(\bar{t}_{\segg{\bar{n}}})\ +\\
&\quad+\sum_{j=2}^{N+1}
\FF^+_{j,1}(t)k^+_j(t)\sk{\EE^+_{1,j}(t)\cdot \F(\bar t_{\segg{\bar n}})-\Aint_{j-1}-
\sum_{m=1}^{j-2}\EE^+_{m+1,j}(t)\cdot \Aint_m}\,,
\end{split}
\end{equation}
where
\begin{equation}\label{Aint}
\Aint_m=(q-q^{-1})^m\ \prod_{a=1}^m (n_a)\ \tSym_{\ \bar t_{\segg{\bar n}}}
\left( \F(\bar t_{\segg{\bar n'}_m})
\cdot \F(\bar t^1_{\segg{n_1-1}})\
\cdot \Zfun_m(\bar t_{\segg{\bar n}})\   \frac{ t^1_{n_1}/t}{1-t^1_{n_1}/t}  \right).
\end{equation}

Calculation of \r{gc8} modulo elements of the ideal $I$ is useful since now to continue ordering
we have to exchange the Gauss coordinates $\EE^+_{m+1,j}$ with total currents. These exchange
relations are based on the formulae given by  Lemma~\ref{5.6e}:
\begin{equation}\label{E-Gau}
\EE^+_{m+1,j}(t)=(q-q^{-1})^{m+2-j}\hS_{m+1}\,\hS_{m+2}\,\cdots\,\hS_{j-2}\sk{\EE^+_{j-1,j}(t)}\,,
\end{equation}
where the screening operators is defined as $\hS_i(B)=E_i[0]\,B-qB\,E_i[0]$, the formulae
\begin{equation}\label{use10}
[\EE^+_{i,i+1}(t),F_{j}(t')]\simK \delta_{i,j}\frac{(q-q^{-1})t'/t}{1-t'/t}\ \psi_i^+(t')\,,
\end{equation}
and
\begin{equation}\label{use11}
[E_{i}[0],F_{j}(t')]\simK \delta_{i,j}(q-q^{-1})\,\psi_i^+(t')\,.
\end{equation}
Here $\psi^+_i(t)=k^+_i(t)\,k^+_{i+1}(t)^{-1}$. Recall that the symbol `$\simK$' means an equality
modulo terms of the ideal $K$ which composed from the elements of $\Uqglnd$ with any mode of the negative
Cartan current $\psi^-_i(t)=k^-_i(t)\,k^-_{i+1}(t)^{-1}$.

The ordering of the first term for $j=2$ in the right hand side of \r{gc8} can be performed as
in the case of the algebra  $\Uqdva$ since the Gauss coordinate $\EE^+_{1,2}(t)$ does not commute only
with the currents $F_1(t^1_\ell)$. This term is equal  to (modulo elements from the ideals $K$ and $J$)
\begin{equation}\label{m=2}
\begin{split}
n_1 (q-q^{-1})\FF^+_{2,1}(t)k^+_2(t)\ &\tSym_{\ \bar t_{\segg{\bar n}}}
\left( {\F(\bar t_{\segg{\bar n'}_1})
\cdot \F(\bar t^1_{\seg{n_1}{1}})}\
\frac{t^1_{1}/t}{1-t^1_{1}/t}\ \psi_1^+(t^1_{1})\right.\ -\\
&\left. -\  {\F(\bar t_{\segg{\bar n'}_1})
\cdot \F(\bar t^1_{\segg{n_1-1}})}\
\frac{t^1_{n_1}/t}{1-t^1_{n_1}/t}\
\prod_{j=1}^{n_2}\frac{q-q^{-1}t^1_{n_1}/t^2_j}{1-t^1_{n_1}/t^2_j}  \right).
\end{split}
\end{equation}

The ordering of the Gauss coordinates $\EE^+_{m+1,j}(t)$ with total currents is more involved.
To perform this ordering we have to use besides
\r{E-Gau}, \r{use10}, \r{use11} also the relation
\begin{equation}\label{use12}
\hS_i\sk{\psi^+_{i+1}(t)}=E_i[0]\,\psi^+_{i+1}(t)-q\,\psi^+_{i+1}(t)\,E_i[0]
=(q-q^{-1})\psi^+_{i+1}(t)\EE^+_i(t)\,.
\end{equation}
Fix the $j$th term of the sum  in the right hand side of \r{gc8} and denote
\begin{equation}\label{nice-id}
\Rfun_j=\EE^+_{1,j}(t)\cdot \F(\bar t_{\segg{\bar n}})-(q-q^{-1})^{j-1}\Aint_{j-1}-
\sum_{m=1}^{j-2} (q-q^{-1})^{m}\EE^+_{m+1,j}(t)\cdot \Aint_m\,.
\end{equation}
We will exchange the Gauss coordinates $\EE^+_{m+1,j}(t)$ with total current.
As above we will calculate modulo elements of the ideals $K$ and $J$.

We will describe this calculation as the sequence of the steps.

\begin{itemize}
\item[{\rm  1.}]
According to \r{E-Gau} the Gauss coordinate $\EE^+_{m+1,j}(t)$ is composed from the
zero mode of the currents $E_i(t)$ with $i=m+1,\ldots,j-2$ and the Gauss coordinate
$\EE^+_{j-1,j}(t)$. It means that this Gauss coordinate will have nontrivial commutation
relations only with the currents $F_a(t^a_s)$ for $a=m+1,\ldots,j-1$.
The move of the Gauss coordinates $\EE^+_{m+1,j}(t)$
through the total current $F_a(t^a_s)$ with $a=m+1,\ldots,j-2$
produces the terms which belong to the ideal $J$. These terms
 appear after commuting zero modes $E_a[0]$ with
the current $F_a(t^a_s)$ and can be neglected since we perform the calculation modulo elements of
the ideal $J$. The nontrivial terms, which we will consider, are created  after exchanging the Gauss coordinate
$\EE^+_{m+1,j}(t)$ with the currents $F_{j-1}(t^{j-1}_s)$.
At the first step the Cartan currents $\psi^+_{j-1}(t^{j-1}_s)$ together with a rational factor
$\frac{(q-q^{-1})t^{j-1}_1/t}{1- t^{j-1}_1/t}$
 appear according to \r{use10}
after an exchange of $\EE^+_{j-1,j}(t)$ and the total current $F_{j-1}(t^{j-1}_s)$ in the
place of the latter current.
This Cartan current
 should be moved to the right  of the product of type $j-1$ currents.
All the terms resulting from the exchange relation of $\EE^+_{j-1,j}(t)$ and  $F_{j-1}(t^{j-1}_s)$
can be presented as the $q$-symmetrization over the type $j-1$ formal variables
 of the single term which appears after commutation of
$\EE^+_{j-1,j}(t)$ with $F_{j-1}(t^{j-1}_1)$. Due to properties \r{sym*} this $q$-symmetrization
is absorbed into `global' $q$-symmetrization entering into definitions of the elements $\Aint_m$
\r{Aint} producing factor $n_{j-1}$ due to definition of the $q$-symmetrization.
The product of the current $\F(\bar t_{\segg{\bar n}})$ in \r{nice-id} can be also
presented as the `global' $q$-symmetrization over the set of the variables $\bar t_{\segg{\bar n}}$.
\item[{\rm 2.}]
The sequence of the screening operators $\hS_{m+1}\cdots \hS_{j-2}$ which enter the
formula for the Gauss coordinate $\EE^+_{m+1,j}(t)$ are applied
 according to equation \r{use12} to replace the Cartan current $\psi^+_{j-1}(t^{j-1}_1)$
on the right of the product of $j-1$-type currents by
the factor $\psi^+_{j-1}(t^{j-1}_1)\EE^+_{m+1,j-1}(t^{j-1}_1)$.
\item[{\rm 3.}]
The Gauss coordinate $\EE^+_{m+1,j-1}(t^{j-1}_1)$ should be moved through the total currents
of type $j-2$ to produce $q$-symmetrization over variables of the same type
and the factor $\psi^+_{j-2}(t^{j-2}_1)\EE^+_{m+1,j-2}(t^{j-2}_1)$ on the right of the product
of this type current accompanied with a rational factor
$\frac{(q-q^{-1})t^{j-2}_1/t^{j-1}_1}{1- t^{j-2}_1/t^{j-1}_{1}}$. Note that the factor $n_{j-2}$
appears also due to definition of the $q$-symmetrization.
\item[{\rm 4.}]
After moving all Gauss coordinates $\EE^+_{m+1,a}(t^{a}_1)$ to the right
we are left with
a product of total currents where all the currents $F_a(t^a_1)$ have
been replaced by
the Cartan currents $\psi_a^+(t^a_1)$ with $a=m+1,\ldots,j-1$ together with the rational
factors $n_a\,\frac{(q-q^{-1})t^{a}_1/t^{a+1}_1}{1- t^{a}_1/t^{a+1}_{1}}$. Then, we have to move
all these Cartan currents to the right of all total currents producing the
rational factors of the Bethe parameters.
\end{itemize}

These calculations result to the following
exchange relations of the Gauss coordinate $\EE^+_{m+1,j}(t)$ and the product of the currents
$\F(\bar t_{\segg{\bar n'}_m})
\cdot \F(\bar t^1_{\segg{n_1-1}})$
\begin{equation*}
\begin{split}
&\EE^+_{m+1,j}(t) \F(\bar t_{\segg{\bar n'}_m})
\cdot \F(\bar t^1_{\segg{n_1-1}}) \simJK \\
&\quad\simJK\ (q-q^{-1})^{j-m-1}\prod_{a=m+1}^{j-1}(n_a)\
\tSym_{\ \bar t_{\segg{\bar n}}}\left(
\F(\bar t^N_{\segg{n_N}})\cdots
\F(\bar t^j_{\segg{n_j}})\ \times\phantom{\prod_a^j}\right.\\
&\qquad\times\ \F(\bar t^{j-1}_{\seg{n_{j-1}}{1}})\cdots \F(\bar t^{m+1}_{\seg{n_{m+1}}{1}})
\F(\bar t^{m}_{\segg{n_{m}-1}})\cdots \F(\bar t^{1}_{\segg{n_{1}-1}})\  \times \\
&\qquad \times \left. \frac{t^{j-1}_1/t}{1-t^{j-1}_1/t}\ \Yfun_m(\bar t_{\segg{\bar n}})\
\prod_{a=m+1}^{j-1}\psi^+_a(t^a_1)\right),
\end{split}
\end{equation*}
where
a rational series $\Yfun_m(\bar t_{\segg{\bar n}})$ is defined by the relation
\begin{equation}\label{Yfun}
\Yfun_m(\bar t_{\segg{\bar n}})=\prod_{a=m+1}^{j-2}\sk{
\frac{t_1^a/t_1^{a+1}}{1-t_1^a/t_1^{a+1}} \prod_{k=2}^{n_a}\frac{q-q^{-1}t_k^a/t_1^{a+1}}{1-t_k^a/t_1^{a+1}}
} \prod_{k=1}^{n_m-1} \frac{q-q^{-1}t_k^m/t_1^{m+1}}{1-t_k^m/t_1^{m+1}}\,.
\end{equation}
Recall that the notation
$\bar t^{a}_{\seg{n_{a}}{1}}$ means the collection of the variables
$\{t_2^{a},\ldots,t^a_{n_a}\}$.

We  get finally the following presentation of the product
$\mathcal{T}_{N+1}(t)\cdot \mathcal{W}_{N+1}(\bar{t}_{\segg{\bar{n}}})$:

\begin{equation}\label{gc88}
\begin{split}
&\mathcal{T}_{N+1}(t)\cdot \mathcal{W}_{N+1}(\bar{t}_{\segg{\bar{n}}})\ \simm\
\mathcal{T}'_{N}(t)\cdot \mathcal{W}'_{N}(\bar{t}_{\segg{\bar{n}'}})\cdot
\F(\bar t_{\segg{n_1}})+\\
&\quad +\mathcal{W}_{N+1}(\bar{t}_{\segg{\bar{n}}})\cdot k_1^+(t)\
\prod_{i=1}^{n_1}\frac{q^{-1}-q t^1_i/t}{1-t^1_i/t}\
\prod_{i=1}^{n_2}\frac{q-q^{-1} t^2_i/t}{1-t^2_i/t}\ +\\
&\quad +\sum_{j=2}^{N+1} (q-q^{-1})^{j-1}\prod_{a=1}^{j-1}(n_a)\
\FF^+_{j,1}(t)k^+_j(t)\cdot \Rfun_j\,,
\end{split}
\end{equation}
where the symbol `$\simm$' means an equality modulo the ideals
$I$, $J$ and $K$ (see definition \ref{def:K}). The elements $\Rfun_j$ have the structure
\begin{equation}\label{nice1}
\begin{split}
\Rfun_j&=\ \tSym_{\ \bar t_{\segg{\bar n}}}
\left ( \F(\bar t^N_{\segg{n_N}})\cdots
\F(\bar t^j_{\segg{n_j}})\F(\bar t^{j-1}_{\seg{n_{j-1}}{1}})\cdots \F(\bar t^{1}_{\seg{n_{1}}{1}})\
 \Yfun_0(\bar t_{\segg{\bar n}})\ \frac{t^{j-1}_1/t}{1-t^{j-1}_1/t}\
 \prod_{a=1}^{j-1}\psi^+_a(t^a_1)\ -  \right.\\
&\qquad-\
\F(\bar t^N_{\segg{n_N}})\cdots
\F(\bar t^j_{\segg{n_j}})\F(\bar t^{j-1}_{\segg{n_{j-1}-1}})\cdots \F(\bar t^{1}_{\segg{n_{1}-1}})
 \cdot \Zfun_{j-1}(\bar t_{\segg{\bar n}})\   \frac{ t^1_{n_1}/t}{1-t^1_{n_1}/t}\ -\\
&\qquad -\ \sum_{m=1}^{j-2} \left(\F(\bar t^N_{\segg{n_N}})\cdots
\F(\bar t^j_{\segg{n_j}})\F(\bar t^{j-1}_{\seg{n_{j-1}}{1}})\cdots \F(\bar t^{m+1}_{\seg{n_{m+1}}{1}})
\F(\bar t^{m}_{\segg{n_{m}-1}})\cdots \F(\bar t^{1}_{\segg{n_{1}-1}})\times\phantom{\prod_{a}^b}\right.\\
&\qquad\qquad \left.\left. \times\  \frac{t^1_{n_1}/t}{1-t^1_{n_1}/t}\ \frac{t^{j-1}_1/t}{1-t^{j-1}_1/t}\
\Zfun_m(\bar t_{\segg{\bar n}})\ \Yfun_m(\bar t_{\segg{\bar n}}) \prod_{a=m+1}^{j-1}\psi^+_a(t^a_1)
 \right)\right).
\end{split}
\end{equation}
The series $\Zfun_m(\bar t_{\segg{\bar n}})$ is defined by the relation \r{Zfun}.

\subsubsection{Proof of Theorem~\ref{main-th}}

The following proposition generalizes
 Proposition~\ref{prop4.2} and serves as a main statement which is proved
 by induction over $N$.
\begin{proposition}\label{4.2gen}
There is a formal series equality in the algebra $\Uqglnd$
\begin{equation}\label{ind1}
\begin{split}
&\mathcal{T}_{N}(t)\cdot \mathcal{W}_{N}(\bar{t}_{\segg{\bar{n}}}) \simm
\mathcal{W}_{N}(\bar{t}_{\segg{\bar{n}}})\cdot
\sk{\sum_{a=1}^{N} k^+_a(t)\prod_{k=1}^{n_{a-1}}\frac{q-q^{-1}t^{a-1}_k/t}{1-t^{a-1}_k/t}
\prod_{k=1}^{n_{a}}\frac{q^{-1}-qt^{a}_k/t}{1-t^{a}_k/t}}\ +\\
&\quad +\tSym_{\ \bar t_{\segg{\bar n}}}\left(\sum_{j=2}^{N}\prod_{a=1}^{j-1}(n_a)\left( (q-q^{-1})^{j-1}\
\FF^+_{j,1}(t)k^+_j(t)\cdot
\F(\bar t^{N-1}_{\segg{n_{N-1}}})\cdots
\F(\bar t^j_{\segg{n_j}})\ \times \phantom{\prod_{a=1}^{j-1}} \right.\right.\\
&\quad\quad\times \left.\F(\bar t^{j-1}_{\segg{n_{j-1}-1}})\cdots \F(\bar t^{1}_{\segg{n_{1}-1}})
 \cdot \Zfun_{j-1}(\bar t_{\segg{\bar n}})\   \frac{ t^{j-1}_{n_{j-1}}/t}{1-t^{j-1}_{n_{j-1}}/t}\right) \times\\
&\quad\quad\times\left.
\sk{\psi^+_1(t^1_{n_1})\prod_{k=1}^{n_1-1}\frac{q^{-1}-qt^1_k/t^1_{n_1}}{q-q^{-1}t^1_k/t^1_{n_1}}
\prod_{k=1}^{n_2}\frac{1-t^2_k/t^1_{n_1}}{q^{-1}-qt^2_k/t^1_{n_1}} -1 }\right)
\end{split}
\end{equation}
 if the set $\{t^i_j\}$ of the Bethe parameters satisfies the set of the universal Bethe equations,
$i=2,\ldots,N-1$, $j=1,\ldots,n_i$:
\begin{equation}\label{u-BE-ind}
\frac{k^+_i(t^i_j)}{k^+_{i+1}(t^i_j)}=
\prod_{m\neq j}^{n_i}\frac{q-q^{-1}t^{i}_m/t^i_j}{q^{-1}-qt^{i}_m/t^i_j}\
\prod_{m=1}^{n_{i-1}}\frac{1-t^{i-1}_m/t^i_j}{q-q^{-1}t^{i-1}_m/t^i_j}\
\prod_{m=1}^{n_{i+1}}\frac{q^{-1}-qt^{i+1}_m/t^i_j}{1-t^{i+1}_m/t^i_j}\,.
\end{equation}
\end{proposition}
\noindent {\it Proof.}\/ We will proof this proposition by induction over $N$ taking as the base
of the induction the statement of Proposition~\ref{prop4.2}.
We assume that the equality \r{ind1} is valid for the algebra $\Uqglnd$ embedded into $\Uqd{N+1}$ by the
relation \r{embed} and prove from the relation \r{gc88}  that \r{ind1} is valid also for the algebra
$\Uqd{N+1}$.

First we rewrite the induction assumption \r{ind1} for the algebra $\Uqglnd$ embedded into $\Uqd{N+1}$ by
 \r{embed}.
It will takes the form
\begin{equation}\label{inemb}
\mathcal{T}'_{N}(t)\cdot \mathcal{W}'_{N}(\bar{t}_{\segg{\bar{n}'}}) \simm
\mathcal{W}'_{N}(\bar{t}_{\segg{\bar{n}'}})\cdot
\sk{\sum_{a=2}^{N+1} k^+_a(t)\prod_{k=1}^{n_{a-1}}\frac{q-q^{-1}t^{a-1}_k/t}{1-t^{a-1}_k/t}
\prod_{k=1}^{n_{a}}\frac{q^{-1}-qt^{a}_k/t}{1-t^{a}_k/t}}\ + \Qfun
\end{equation}
where in `unwanted' terms
\begin{equation}\label{ind1emb}
\begin{split}
&\Qfun=\tSym_{\ \bar t_{\segg{\bar n}}}\left(\sum_{j=3}^{N+1}\left( (q-q^{-1})^{j-1}\prod_{a=2}^{j-1}(n_a)\
\FF^+_{j,2}(t)k^+_j(t)\cdot
\F(\bar t^{N}_{\segg{n_{N}}})\cdots
\F(\bar t^j_{\segg{n_j}})\ \times \right.\right.\\
&\quad\quad\times \left.\F(\bar t^{j-1}_{\segg{n_{j-1}-1}})\cdots \F(\bar t^{2}_{\segg{n_{2}-1}})
 \cdot \Zfun_{j-1}(\bar t_{\segg{\bar n'}})\   \frac{ t^{j-1}_{n_{j-1}}/t}{1-t^{j-1}_{n_{j-1}}/t}\right) \times\\
&\quad\quad\times\left.
\sk{\psi^+_2(t^2_{n_2})\prod_{k=1}^{n_2-1}\frac{q^{-1}-qt^2_k/t^2_{n_2}}{q-q^{-1}t^2_k/t^2_{n_2}}
\prod_{k=1}^{n_3}\frac{1-t^3_k/t^2_{n_2}}{q^{-1}-qt^3_k/t^2_{n_2}} -1 }\right)
\end{split}
\end{equation}
parameters $\bar t^a_{\bar n_a}$ with $a=3,\ldots,N$ satisfy the universal Bethe equations \r{u-BE-ind}
for $i=3,\ldots,N$ while
parameters $\bar t^2_{\bar n_2}$ are free.

We substitute \r{inemb} into \r{gc88}. First term in the right hand side of \r{inemb} together with the
second term in \r{gc88} produces `wanted' terms
\begin{equation*}
\mathcal{W}_{N+1}(\bar{t}_{\segg{\bar{n}}})\cdot
\sk{\sum_{a=1}^{N+1} k^+_a(t)\prod_{k=1}^{n_{a-1}}\frac{q-q^{-1}t^{a-1}_k/t}{1-t^{a-1}_k/t}
\prod_{k=1}^{n_{a}}\frac{q^{-1}-qt^{a}_k/t}{1-t^{a}_k/t}}.
\end{equation*}
The terms in the right hand side of \r{inemb} which belongs to the ideal $J$ will be again in the same ideal
after multiplication to the right  of $\mathcal{T}'_{N}(t)\cdot \mathcal{W}'_{N}(\bar{t}_{\segg{\bar{n}'}})$ by the
product of the first type currents $\F(\bar t_{\segg{n_1}})$ (see \r{gc88}).
This is because the currents $E_a(t)$ commute with the current $F_1(t')$ for $a=2,\ldots,N$.

Fix parameters $\bar t^2_{\bar n_2}$ from the condition that the
ordered product $:\Qfun\cdot \F(\bar t^1_{\segg{n_1}}):$ of the unwanted terms and the currents of the
first type $F_1(t^1_k)$ vanishes. This results into Bethe equations
\begin{equation}\label{BEmore}
\frac{k^+_2(t^2_k)}{k^+_{3}(t^2_k)}=
\prod_{m\neq k}^{n_2}\frac{q-q^{-1}t^{2}_m/t^2_k}{q^{-1}-qt^{2}_m/t^2_k}\
\prod_{m=1}^{n_{1}}\frac{1-t^{1}_m/t^2_k}{q-q^{-1}t^{1}_m/t^2_k}\
\prod_{m=1}^{n_{3}}\frac{q^{-1}-qt^{3}_m/t^2_k}{1-t^{3}_m/t^2_k}
\end{equation}
for the set of parameters $\bar t^2_{\bar n_2}$.

Now we will examine the structure of the terms $\Rfun_j$ given by \r{nice1} by the conditions that
parameters $\bar t^a_{\bar n_a}$ with $a=2,\ldots,N$ are bounded by the universal Bethe equations
\r{u-BE-ind} and \r{BEmore}.
We replace in $\Rfun_j$ the Cartan currents
$\psi_a^+(t^a_1)$, $a=2,\ldots,j-1$  by the right hand sides of the universal Bethe equations.
Each Bethe equation introduces the factor
$\prod_{\ell=2}^{n_a}\frac{qt^a_1-q^{-1}t^a_\ell}{q^{-1}t^a_1-qt^a_\ell}$ under $q$-symmetrization.
This allows to use the following property of the $q$-symmetrization, which is a consequence of
\r{sym-last} and \r{sym-first}:
\begin{equation}\label{qsy-pr}
\tSym_{\ \bar t}\sk{G(t_1,t_2,\ldots,t_n)\prod_{\ell=2}^n
\frac{q^{-1}-qt_1/t_\ell}{q-q^{-1}t_1/t_\ell}}=
\tSym_{\ \bar t}\sk{G(t_n,t_1,\ldots,t_{n-1})}
\end{equation}
for arbitrary formal series $G$ of the formal variables $(t_1,\ldots,t_n)$.
The  variables $\{t^a_1,t^a_2,\ldots,t^a_{n_a}\}$ are replaced by the permuted
sets $\{t^a_{n_a},t^a_1,\ldots,t^a_{n_a-1}\}$ for $a=m+1,\ldots,j-1$. Using an
identity
\begin{equation*}
\begin{split}
&
\frac{1}{t-t^{j-1}_{n_{j-1}}}\ \prod_{a=1}^{j-2} \frac{1}{t^{a+1}_{n_{a+1}}-t^{a}_{n_{a}}}
-\frac{1}{t-t^{1}_{n_{1}}}\ \prod_{a=1}^{j-2} \frac{1}{t^{a+1}_{n_{a+1}}-t^{a}_{n_{a}}}
\\
&\quad-\frac{1}{t-t^{j-1}_{n_{j-1}}}\ \frac{1}{t-t^{1}_{n_{1}}}\
\sum_{m=1}^{j-2}\prod_{a=m+1}^{j-2} \frac{1}{t^{a+1}_{n_{a+1}}-t^{a}_{n_{a}}}
\prod_{a=1}^{m-1} \frac{1}{t^{a+1}_{n_{a+1}}-t^{a}_{n_{a}}}=0\,.
\end{split}
\end{equation*}
and the explicit forms of the series $\Zfun_m$ and $\Yfun_m$ we observe that
the element $\Rfun_j$ has the structure
\begin{equation*}
\begin{split}
&\Rfun_j=\F(\bar t^{N}_{\segg{n_{N}}})\cdots
\F(\bar t^j_{\segg{n_j}})\cdot\F(\bar t^{j-1}_{\segg{n_{j-1}-1}})\cdots \F(\bar t^{1}_{\segg{n_{1}-1}})
 \cdot \Zfun_{j-1}(\bar t_{\segg{\bar n}})\   \frac{ t^{j-1}_{n_{j-1}}/t}{1-t^{j-1}_{n_{j-1}}/t}\times\\
&\quad\times
\sk{\psi^+_1(t^1_{n_1})\prod_{k=1}^{n_1-1}\frac{q^{-1}-qt^1_k/t^1_{n_1}}{q-q^{-1}t^1_k/t^1_{n_1}}
\prod_{k=1}^{n_2}\frac{1-t^2_k/t^1_{n_1}}{q^{-1}-qt^2_k/t^1_{n_1}} -1 }
\end{split}
\end{equation*}
and this proves the statement of the Proposition.\qed

The statement  of Theorem~\ref{main-th} follows from Proposition~\ref{4.2gen} if we impose
one more universal Bethe equation for the set of the parameters $\{\bar t^1_{\segg{n_1}}\}$
\begin{equation}\label{BE-last}
\frac{k^+_1(t^1_j)}{k^+_{2}(t^1_j)}=
\prod_{m\neq j}^{n_1}\frac{q-q^{-1}t^{1}_m/t^1_j}{q^{-1}-qt^{1}_m/t^1_j}\
\prod_{m=1}^{n_{2}}\frac{q^{-1}-qt^{2}_m/t^1_j}{1-t^{2}_m/t^1_j}\,.
\end{equation}
Since the left hand side and first term in the right hand side of \r{ind1} belong to the standard
Borel subalgebra $\Uqbp$ in $\Uqglnd$ the equality between them is valid modulo elements of the ideal $J$
only if all the parameters satisfies the universal Bethe equations \r{BE-last} and
\r{u-BE-ind}. This finishes the proof of Theorem~\ref{main-th}.\qed

\section*{Acknowledgement}

This work was partially done when the third
author (S.P.) visited Laboratoire d'An\-necy-Le-Vieux de Physique Th\'eorique in 2006 and 2007.
These visits were possible due to the financial support of
the CNRS-Russia exchange program on mathematical physics.
He thanks LAPTH for the hospitality and stimulating scientific atmosphere.
Work  by second (S.K) and third (S.P.) authors was supported in part by RFBR
grant 08-01-00667 and grant for support of scientific schools NSh-3036.2008.2.
Third author (S.P.) was also supported in part by RFBR-CNRS grant 07-02-92166-CNRS.

\end{document}